\documentclass[aap,MSNbibl,nameyear,seceqn,dvips]{arximspdf}
\usepackage{graphicx}
\doi{10.1214/13-AAP934} 
\volume{24}
\issue{2}
\pubyear{2014}
\firstpage{721}
\lastpage{759}

\makeatletter
\def\Bbb{\mathbb}
\def\boldsymbol{\bolds}

\def\cal{\mathcal}

\newtheorem{theorem}{Theorem}[section]
\newtheorem{lemma}[theorem]{Lemma}
\newtheorem{condition}[theorem]{Condition}

\newproclaim{remark}[theorem]{Remark}

\newproclaim{example}[theorem]{Example}
\newproclaim{definition}[theorem]{Definition}

\newcommand{\I}{\mathcal{I}}
\newcommand{\K}{\mathcal{K}}

\def\hat{\widehat}
\def\tilde{\widetilde}
\def\bar{\overline}
\def\del{\partial}
\makeatother

\begin{document}
\begin{frontmatter}

\title{Central limit theorems and diffusion approximations for
multiscale Markov~chain~models\thanksref{T1}}
\runtitle{CLT for Multi-scale Markov chains}
\thankstext{T1}{Supported in part by NSF FRG Grant DMS 05-53687.}

\begin{aug}
\author[A]{\fnms{Hye-Won} \snm{Kang}\ead[label=e1]{kang235@mbi.osu.edu}},
\author[B]{\fnms{Thomas G.} \snm{Kurtz}\thanksref{t2}\ead[label=e2]{kurtz@math.wisc.edu}}
\and
\author[C]{\fnms{Lea} \snm{Popovic}\corref{}\thanksref{t3}\ead[label=e3]{lpopovic@mathstat.concordia.ca}}
\thankstext{t2}{Supported in part by NSF Grant DMS 11-06424.}
\thankstext{t3}{Supported by NSERC Discovery grant.}
\runauthor{H.-W. Kang, T. G. Kurtz and L. Popovic}
\affiliation{Ohio State University, University of Wisconsin
and Concordia University}
\address[A]{H.-W. Kang\\
Mathematical Biosciences Institute\\
Ohio State University\\
1735 Neil Ave.\\
Columbus, Ohio 43210\\
USA\\
\printead{e1}} 
\address[B]{T. G. Kurtz\\
Departments of Mathematics and Statistics\\
University of Wisconsin\\
480 Lincoln Dr.\\
Madison, Wisconsin 53706\\
USA\\
\printead{e2}}
\address[C]{L. Popovic\\
Department of Mathematics and Statistics\\
Concordia University\\
1450 de Maisonneuve Blvd. West\\
Montreal, Quebec H3G1M8\\
Canada\\
\printead{e3}}
\end{aug}

\received{\smonth{8} \syear{2012}}
\revised{\smonth{3} \syear{2013}}

%
\begin{abstract}
Ordinary differential equations obtained as limits of Markov processes
appear in many settings. They may arise by scaling large systems, or by
averaging rapidly fluctuating systems, or in systems involving multiple
time-scales, by a combination of the two. Motivated by models with
multiple time-scales arising in systems biology, we present a general
approach to proving a central limit theorem capturing the fluctuations
of the original model around the deterministic limit. The central limit
theorem provides a method for deriving an appropriate diffusion
(Langevin) approximation.
\end{abstract}

%

\end{frontmatter}

\section{Introduction}\label{intro}
There are two classical kinds of Gaussian limit theorems
associated with continuous time Markov chains as well
as more general Markov processes. The first of these
considers a sequence $\{X^N\}$ of Markov chains that
converges to a deterministic function $X$ and gives a limit
for the rescaled deviations $U^N=r_N(X^N-X)$; see, for
example, \citet{Kur71,Ku77,vanK61}. The second
considers an ergodic Markov process $Y$ with stationary
distribution $\pi$ and gives a limit
for
\[
Z^N(t)=\frac{1}{\sqrt{N}}\int_0^{Nt}h
\bigl(Y(s)\bigr)\,ds=\sqrt{N}\int_0^t h\bigl(Y(Ns)
\bigr)\,ds
\]
for $h$ satisfying $\int h\,d\pi=0$; see, for example, \citet{Bha82} for
a general result of
this type.

There are many proofs for theorems like these. In
particular, results of both types can be proved using the
martingale central limit theorem (Theorem \ref{mgclt}).
For example, in the first case, there is typically a
sequence of functions $F^N$ such that
\[
M^N(t)=X^N(t)-X^N(0)-\int
_0^tF^N\bigl(X^N(s)
\bigr)\,ds
\]
is a martingale, $F^N\rightarrow F$,
$\dot{X}=F(X)$ and $F^N(x^N)-F(x)\approx\nabla F(x)(x^N-x)$, for $
x^N$
converging to $x$. If the martingale central limit theorem
gives $r_NM^N\Rightarrow M$ and $U^N(0)\Rightarrow U(0)$, then (ignoring
technicalities) $U^N$ should converge to the solution of
%
\begin{equation}
U(t)=U(0)+M(t)+\int_0^t\nabla F\bigl(X(s)
\bigr)U(s)\,ds.\label{ou}
\end{equation}

In the second case, the assumption that $\int h\,d\pi=0$ suggests
that there should be a solution of the Poisson equation
$Af=-h$, where $A$ is the generator for $Y$, and then
\begin{eqnarray*}
Z^N(t)&=&\frac{1}{\sqrt{N}}\biggl(f\bigl(Y(Nt)\bigr)-f\bigl(Y(0)\bigr)-
\int_0^{Nt}Af\bigl(Y(s)\bigr)\,d s\biggr)\\
&&{}-
\frac{1}{\sqrt{N}}\bigl(f\bigl(Y(Nt)\bigr)-f\bigl(Y(0)\bigr)\bigr).
\end{eqnarray*}
The first term on the right is a martingale and the
second should go to zero, so if the
martingale central
limit theorem applies to the first, then $Z^N$ should
converge.

This paper addresses situations of the first type
[$V_0^N\Rightarrow V_0$ for a deterministic~$V_0$, and we want to verify
convergence of $U^N=r_N(V_0^N-V_0)$] in which both approaches
are required. Specifically, the function $F^N$ giving the
martingale, $M^N$, depends not
only on $V_0^N$ but also on another process $V_1^N$ [think
$V_1^N(t)=V_1(Nt)$], so
\[
M^{N,1}(t)=V_0^N(t)-V_0^N(0)-
\int_0^tF^N\bigl(V_0^N(s),V_1^N(s)
\bigr)\,ds
\]
is a martingale, $F^N$ ``averages'' to $F$ in the sense that
\[
\int_0^t\bigl(F^N
\bigl(V_0^N(s),V_1^N(s)\bigr)-F
\bigl(V_0^N(s)\bigr)\bigr)\,ds\rightarrow0,
\]
$(V_0^N,V_1^N)$ is Markov with generator $A_N$ and there exist
$H_N$
such that $A_NH_N\approx(F^N-F)$. [Note that $H_N$ will be a
vector of functions in the domain of $A_N$, ${\cal D}(A_N)$.]
Assuming that
\begin{eqnarray*}
M^{N,2}(t)&=&H_N\bigl(V_0^N(t),V_1^N(t)
\bigr)-H_N\bigl(V_0^N(0),V_1^N(0)
\bigr)\\
&&{}-\int_0^
tA_NH_N
\bigl(V_0^N(s),V_1^N(s)\bigr)\,ds
\end{eqnarray*}
is a martingale, and
again ignoring all the
technicalities, we have
%
\begin{eqnarray}\label{beq}
&&r_N\bigl(V_0^N(t)-V_0(t)
\bigr)\nonumber
\\
&&\qquad =r_N\bigl(V_0^N(0)-V_0(0)
\bigr)+r_NM^{N,1}(t)-r_NM^{N,2}(t)
\nonumber
\\
&&\qquad\quad{} +\int_0^tr_N\bigl(F
\bigl(V_0^N(s)\bigr)-F\bigl(V_0(s)\bigr)
\bigr)\,ds
\\
&&\qquad\quad{} +r_N\bigl(H_N\bigl(V_0^N(t),V_1^N(t)
\bigr)-H_N\bigl(V_0^N(0),V_1^N(0)
\bigr)\bigr)
\nonumber\\
&&\qquad\quad{} +r_N\int_0^t
\bigl(F^N\bigl(V_0^N(s),V_1^N(s)
\bigr)-F\bigl(V_0^N(s)\bigr)\nonumber\\
&&\hspace*{75pt}\qquad\quad{}{}-A_NH_N
\bigl(V_
0^N(s),V_1^N(s)\bigr)
\bigr)\,ds.
\nonumber
\end{eqnarray}
If the last two terms on the right go to zero, the
martingale terms converge,
\[
r_NM^{N,1}-r_NM^{N,2}\Rightarrow M
\]
and $F$ is
smooth, then we again should have $U^N\Rightarrow U$ satisfying
(\ref{ou}).

The work to be done to obtain theorems of this type is
now clear. We need to identify $F^N$ and $F$, find an
approximate solution to the Poisson equation
$A_NH_N\approx F^N-F$, verify that the martingales satisfy the
conditions of the martingale central limit theorem, and
verify that the error terms [the last two terms in
(\ref{beq})] converge to zero. We will make this analysis
more specific in stages. We are essentially considering
situations in which the process $V_1^N$ is evolving on a
faster time scale than $V^N_0$ and ``averages out'' to give the
convergence of $V^N_0$ to $V_0$. But $V_1^N$ itself may evolve on
more than one time scale. In the first stage of our
development, we will replace $V_1^N$ by $(V_1^N,V_2^N)$ with $V_1^
N$
and $V_2^N$ evolving on different (fast) time scales. Once
the analysis for two fast time scales is carried out, the
extension of the general results to more than two fast
time scales should be clear. In the second stage, we
consider multiply scaled, continuous-time Markov chains
of a type that arises naturally in models of chemical
reaction networks. For these models, many of the
conditions simplify, but the notation becomes more
complex.

\textit{Outline}:
In Section~\ref{sectclt} we state and prove the functional central
limit Theorem~\ref{genclt}, and specify a sequence of Conditions~\ref
{scparam}--\ref{tight} that need to be verified for it to apply. In
Section~\ref{sectdiff} we additionally give a diffusion approximation
implied by Theorem~\ref{genclt}. Our aim is to apply these results to
Markov chain models for chemical reactions. In Section~\ref{sectmcmcr}
we identify specific aspects of the multi-scale behavior of a reaction
network that one needs in order to apply Theorem~\ref{genclt} to the
chemical species with a deterministic limit on the slowest time scale.
Section~\ref{sectexamp} provides several examples of chemical networks
(the first two evolving on two, the last one on three time-scales), and
shows how to verify the conditions and obtain a diffusion approximation.

\section{A central limit theorem for a system with
deterministic limit and three time scales}\label{sectclt}

We identify a set of conditions on a three time-scale process
$V^N=(V_0^N, V_1^N, V_2^N)$ that guarantee $U^N=r_N(V_0^N-V_0)$
converges to a diffusion. As suggested earlier, we write $U^N$ in the form
%
\begin{eqnarray}\label{uexp}
U^N(t)&=&U^N(0)+r_N\bigl(M^{N,1}(t)-M^{N,2}(t)
\bigr)\nonumber
\\
&&{} +r_N\int_0^t\bigl(\bar{F}
\bigl(V_0^N(s)\bigr)-\bar{F}\bigl(V_0(s)\bigr)
\bigr)\,ds
\\
&&{} +r_N\int_0^t
\bigl(F^N\bigl(V^N(s)\bigr)-F\bigl(V^N(s)
\bigr)\bigr)\,ds
\nonumber\\
&& {}+r_N\int_0^t\bigl(F
\bigl(V^N(s)\bigr)-\bar{F}\bigl(V^N_0(s)
\bigr)-A_NH_N\bigl(V^N(s)\bigr) \bigr)\,ds
\nonumber
\\
&&{} +r_N\bigl(H_N\bigl(V^N(t)
\bigr)-H_N\bigl(V^N(0)\bigr)\bigr),
\nonumber
\end{eqnarray}
where $M^{N,1}, M^{N,2}$ are martingales, $V_0$ is the deterministic
limit of the process $V_0^N$, $\bar F$ is its infinitesimal drift and
$H_N$ is an approximate solution to a Poisson equation.
Our conditions insure each individual term has a well behaved limit.

We assume that $V_i^N$ takes values in ${\Bbb E}_i^N\subset{\Bbb R}^{
d_i}$, $i=0,1,2$,
and that ${\Bbb E}_i^N$ converges in the sense that there exists
${\Bbb E}_i\subset{\Bbb R}^{d_i}$ such that ${\Bbb E}_i^N\subset
{\Bbb E}_i$ and for each compact
$K\subset{\Bbb R}^{d_i}$,
\[
\lim_{N\rightarrow\infty}\sup_{x\in{\Bbb E}_i\cap K}\inf
_{y\in
{\Bbb E}_i^N}|x-y|=0.
\]

We will refer to $A_N$ as the ``generator'' for the process
$V^N=(V_0^N,V_1^N,V_2^N)$, but all we require is that $A_N$ is a
linear operator on some space ${\cal D}(A_N)$ of measurable functions
on ${\Bbb E}^N\equiv{\Bbb E}_0^N\times{\Bbb E}_1^N\times{\Bbb E}_
2^N$ and that for $h\in{\cal D}(A_N)$,
\[
h\bigl(V^N(t)\bigr)-h\bigl(V^N(0)\bigr)-\int
_0^tA_Nh\bigl(V^N(s)
\bigr)\,ds
\]
is a local martingale.

We first identify the \textit{time scales of the process} $V^N$ with two
sequences of
positive numbers $\{r_{1,N}\}$, $\{r_{2,N}\}$, and introduce a
sequence of scaling parameters $\{r_N\}$ for $U^N$ with the following
properties.

\begin{condition}[(Scaling
parameters)]\label{scparam} The scaling parameters
$r_N\rightarrow\infty$ and $\{r_{1,N}\}$, $\{r_{2,N}\}$ are sequences
of positive
numbers satisfying
%
\begin{eqnarray}\label{rates}
\lim_{N\rightarrow\infty}
\frac{r_N}{r_{1,N}}&=&0,
\nonumber
\\[-8pt]
\\[-8pt]
\nonumber
\lim_{N\rightarrow\infty}\frac{r_{1,N}}{r_{2,N}}&=&0.
\end{eqnarray}
\end{condition}

We next identify the ``generators'' for the \textit{effective dynamics of
$V_0^N, V_1^N$ and $V_2^N$} on time scales $t, t{r_{1,N}}$, and
$t{r_{2,N}}$, respectively. $L_0$, $L_1$, $L_2$ will be linear operators
defined on sufficiently large domains, ${\cal D}(L_0)\subset M({\Bbb E}_
0)$,
${\cal D}(L_1)\subset M({\Bbb E}_0\times{\Bbb E}_1)$ and ${\cal D}
(L_2)\subset M({\Bbb E}_0\times{\Bbb E}_1\times{\Bbb E}_2)$, and taking
values in $M({\Bbb E}_0\times{\Bbb E}_1\times{\Bbb E}_2)$.
The requirements that
determine what is meant by
``sufficiently large'' will become clear, but we will
assume that the domains contain all $C^{\infty}$ functions having
compact support in the appropriate space.
We will use the notation ${\Bbb E}={\Bbb E}_0\times{\Bbb E}_1\times
{\Bbb E}_2$.

\begin{condition}[(Multiscale convergence)]\label{multconv}
$\!\!\!$For each compact
$K\subset{\Bbb R}^{d_0+d_1+d_2}$,
\begin{eqnarray*}
\lim_{N\rightarrow\infty}\sup_{v\in K\cap{\Bbb E}^N}\bigl|A_Nh(v)-L_
0h(v)\bigr|&=&0,\qquad
h\in{\cal D}(L_0),
\\
\lim_{N\rightarrow\infty}\sup_{v\in K\cap{\Bbb E}^N}\biggl|\frac{1}{r_{
1,N}}A_Nh(v)-L_1h(v)\biggr|&=&0,\qquad
h\in{\cal D}(L_1)
\end{eqnarray*}
and
\[
\lim_{N\rightarrow\infty}\sup_{v\in K\cap{\Bbb E}^N}\biggl|\frac{1}{r_{
2,N}}A_Nh(v)-L_2h(v)\biggr|=0,
h\in{\cal D}(L_2).
\]
\end{condition}

\begin{remark}
Similar conditions are considered in \citet{EN80}. See also
\citet{EK86}, Section~1.7.
There may be only two time-scales, in which case
$d_2=0$, $L_2h=0$ and ${\Bbb E}={\Bbb E}_0\times{\Bbb E}_1$ (equivalently,
${\Bbb E}_2$
consists of a single point) in what follows.
\end{remark}

The next condition ensures the \textit{uniqueness of the
conditional equilibrium distributions of the fast
components $V_2^N$ and $V_1^N$}, whose ``generators'' are $L_2$ and~$L_1$.

\begin{condition}[(Averaging condition)]\label{qstat}
For each $(v_0,v_1)\in{\Bbb E}_0\times{\Bbb E}_1$, there exists a unique
$\mu_{v_0,v_1}\in{\cal P}({\Bbb E}_2)$ such that $\int L_2h(v_0,
v_1,v_2)\mu_{v_0,v_1}(dv_2)=0$ for
\mbox{every} $h\in{\cal D}(L_2)\cap B({\Bbb E})$. For each $v_0\in{\Bbb E}_
0$, there exists a unique
$\mu_{v_0}\in{\cal P}({\Bbb E}_1)$ such that $\int L_1h(v_0,v_1,
v_2)\mu_{v_0,v_1}(dv_2)\mu_{v_0}(dv_1)=0$
for every $h\in{\cal D}(L_1)\cap\break  B({\Bbb E}_0\times{\Bbb E}_1)$.
\end{condition}

With this condition in mind, we define
\[
\bar{L}_1h(v_0,v_1)=\int
L_1h(v_0,v_1,v_2)
\mu_{v_0,v_1}(dv_2).
\]

Our first convergence condition insures that the \textit{slow component
$V_0^N$ has a deterministic limit}.
Essentially it implies that its ``generator'' $L_0h=F\cdot\nabla h$,
for $h\in
C^{\infty}_c({\Bbb E}_0)$.
It also identifies the \textit{intrinsic fluctuations of the slow
component} via a martingale $M^{N,1}$.
For an ${\Bbb R}^{d_0}$-valued process $Y$, we use $[Y]_t$ to denote
the matrix of
covariations $[Y_i,Y_j]_t$.

\begin{condition}[(First convergence condition)]\label{cnd1}
There exist $F^N\in M({\Bbb E}^N,\break  {\Bbb R}^{d_0})$ and $F,G_0\in
C({\Bbb E},{\Bbb R}^{d_0})$ such that
%
\begin{equation}
M^{N,1}(t)=V_0^N(t)-V_0^N(0)-
\int_0^tF^N\bigl(V^N(s)
\bigr)\,d s\label{mg1}
\end{equation}
is a local martingale,
$[V^N_0]_t\Rightarrow0$,
and for each compact $K\subset{\Bbb E}$,
%
\begin{equation}
\lim_{N\rightarrow\infty}\sup_{v\in K\cap{\Bbb E}^
N}\bigl|r_N
\bigl(F^N(v)-F(v)\bigr)-G_0(v)\bigr|=0.\label{err1}
\end{equation}
\end{condition}

We next turn to the relevant \textit{Poisson equations} based on the
conditional equilibrium distributions of the fast components and the
limiting drift of the slow component.
Suppose that there exist $h_1\in{\cal D}(L_1)^{d_0}$ and $h_2,h_
3\in{\cal D}(L_2)^{d_0}$ such that
%
\begin{eqnarray}\label{hnsol}
\bar{L}_1h_1(v_0,v_1)&=&\int
F(v_0,v_1,v_2)\mu_{v_0,v_1}(dv_2)\nonumber\\
&&{}-
\int\!\!\int F(v_0,v_1,v_2)\mu_{v_0,v_1}(dv_2)
\mu_{v_0}(dv_1),
\nonumber
\\[-8pt]
\\[-8pt]
\nonumber
L_2h_2(v_0,v_1,v_2)&=&F(v_0,v_1,v_2)-
\int F(v_0,v_1,v_2)\mu_{v_0,
v_1}(dv_2),
\\
L_2h_3(v_0,v_1,v_2)&=&
\bar {L}_1h_1(v_0,v_1)-L_1h_1(v_0,v_1,v_2).
\nonumber
\end{eqnarray}
Define
\begin{eqnarray*}
\bar{F}_1(v_0,v_1)&=&\int
F(v_0,v_1,v_2)\mu_{v_0,v_1}(dv_2),
\\
\bar{F}(v_0)&=&\int\!\!\int F(v_0,v_1,v_2)
\mu_{v_0,v_1}(dv_2)\mu_{v_0} (dv_1)
\end{eqnarray*}
and
%
\begin{equation}
H_N=\frac{1}{r_{1,N}}h_1+\frac{1}{r_{2,N}}(h_2+h_3)
.\label{hndef}
\end{equation}
Note that for $H_N$ of this form
\[
A_NH_N\approx L_1h_1+L_2(h_2+h_3)=F-
\bar{F}.
\]
In what follows, $H_N$ does not have to be given by
(\ref{hndef}). That form simply suggests the possibility
of finding $H_N$ with the desired properties. Specifically,
we assume the existence of $H_N\in{\cal D}(A_N)$ satisfying the
following convergence condition.

\begin{condition}[(Second convergence condition)]\label{cnd7}
Assume that there exists $G_1\in C({\Bbb E},{\Bbb R}^{d_0})$ such that for
each compact $K\subset{\Bbb E}$,
%
\begin{equation}
\lim_{N\rightarrow\infty}\sup_{v\in K\cap{\Bbb E}^
N}\bigl|r_N
\bigl(F(v)-\bar{F}(v_0)-A_NH_N(v)
\bigr)-G_1(v)\bigr|=0.\label{err2}
\end{equation}
\end{condition}

\begin{remark} The critical
requirements for $H_N$ are (\ref{err2}), (\ref{hneg}) and (\ref{cov1}).
In fact, because of the possibility of large fluctuations
by $V^N_1$ and $V^N_2$, even if $h_1$, $h_2$ and $h_3$ satisfying Condition
\ref{hnsol} can be found,
it may be necessary to define $H_N$ using
a sequence of truncations of $h_1$, $h_2$ and
$h_3$.
\end{remark}

This now identifies the \textit{fluctuations of the slow component due to
convergence of the fast components} to their conditional equilibrium
distributions via a martingale $M^{N,2}$.
For $V_0(0)\in{\Bbb R}^{d_0}$, let $V_0$ satisfy
%
\begin{equation}
V_0(t)=V_0(0)+\int_0^t
\bar{F}\bigl(V_0(s)\bigr)\,ds,\label{limdeq}
\end{equation}
and define
\[
M^{N,2}(t)=H_N\bigl(V^N(t)\bigr)-H_N
\bigl(V^N(0)\bigr)-\int_0^tA_NH_N
\bigl(V^N(s)\bigr)\,ds.
\]

The following condition
is then needed for \textit{application of the martingale central limit
theorem}, Theorem \ref{mgclt}, to the process $M^{N,1}-M^{N,2}$,
composed of $M^{N,2}$ above and $M^{N,1}$ from (\ref{mg1}).
Essentially it says that the jumps of both the slow component and
solutions to the Poisson equations are appropriately small, and that
the quadratic variation of $M^{N,1}-M^{N,2}$ converges.

\begin{condition}[(Converegence of covariation)]\label{hmgcnv}
There exists $G\in C({\Bbb E},\break  {\Bbb M}^{d_0\times d_0})$ such that for
each $
t>0$,
%
\begin{eqnarray}
\lim_{N\rightarrow\infty}E\Bigl[\sup_{s\leq t}r_N\bigl|V_0^
N(s)-V_0^N(s-)\bigr|
\Bigr]&=&0,\label{jmp1}
\\
\sup_{s\leq t}r_NH_N\bigl(V^N(s)
\bigr)&\Rightarrow&0\label{hneg}
\end{eqnarray}
and
%
\begin{equation}
(r_N)^2\bigl[V_0^N-H_N
\circ V^N\bigr]_t-\int_0^tG
\bigl(V^N(s)\bigr)\,ds \Rightarrow0.\label{cov1}
\end{equation}
\end{condition}

We can now \textit{account for all the terms} in the expansion (\ref
{uexp}) of $U^N=r_N(V_0^N-V_0)$,
\begin{eqnarray*}
U^N(t)&=&U^N(0)+r_N\bigl(M^{N,1}(t)-M^{N,2}(t)
\bigr)
\\
&&{} +r_N\int_0^t\bigl(\bar{F}
\bigl(V_0^N(s)\bigr)-\bar{F}\bigl(V_0(s)\bigr)
\bigr)\,ds
\\
&&{} +r_N\int_0^t
\bigl(F^N\bigl(V^N(s)\bigr)-F\bigl(V^N(s)
\bigr)\bigr)\,ds
\\
&&{} +r_N\int_0^t\bigl(F
\bigl(V^N(s)\bigr)-\bar{F}\bigl(V^N_0(s)
\bigr)-A_NH_N\bigl(V^N(s)\bigr) \bigr)\,ds
\\
&& {} +r_N\bigl(H_N\bigl(V^N(t)
\bigr)-H_N\bigl(V^N(0)\bigr)\bigr).
\end{eqnarray*}
Conditions \ref{scparam}--
\ref{hmgcnv} insure that all the terms will have a limit as $N\to
\infty$.
The limit of the second term on the right is guaranteed by (\ref
{jmp1}), (\ref{hneg}) and (\ref{cov1}) in Condition~\ref{hmgcnv}.
Assuming that $\bar{F}$ is smooth,
the third term on the right is asymptotic to $\int_0^t\nabla\bar
{F}(V_0(s))\cdot
U^N(s)\,ds$.
The fourth term is controlled by (\ref{err1}) of Condition~\ref{cnd1},
the fifth by (\ref{err2}) of Condition~\ref{cnd7}. The order of time
scale parameters
(\ref{rates}) from Condition~\ref{scparam} and the form of $H_N$ in
(\ref{hndef}) suggests that the sixth term goes to zero,
but we will explicitly assume that in the statement of the theorem.

Finally, we now only need a condition to \textit{ensure relative
compactness} of the sequence. If $\Bbb E$ is unbounded, let
$\psi\dvtx {\Bbb E}\rightarrow[1,\infty)$ be locally bounded and satisfy
$\lim_{v\rightarrow\infty}\psi(v)=\infty$, \emph{or let
$\psi\dvtx {\Bbb E}\rightarrow[1,\infty)$ be such that $\forall M<\infty
$, $\{v\in\Bbb E\dvtx \psi(v)\le M\}$ is relatively compact in $\Bbb E$,}
and let $D_{\psi}$ denote the collection of
continuous
functions $f$ satisfying
\[
\sup_{v\in{\Bbb E}}\frac{|f(v)|}{\psi(v)}<\infty,\qquad \lim_{
k\rightarrow\infty}
\sup_{v\in{\Bbb E},|v|>k}\frac{|f(v)|}{\psi
(v)}=0.
\]

For sequences of space--time random measures, the
notion of convergence that we will use is that discussed
in \citet{Ku92}.

\begin{lemma}\label{relcomp}
Let $V^N$ be a sequence of ${\Bbb E}$-valued processes, and define
the occupation measure
%
\begin{equation}
\Gamma_N\bigl(D\times[0,t]\bigr)=\int_0^t{
\bf1}_D\bigl(V^N(s)\bigr) \,ds.\label{occdef}
\end{equation}
Suppose that for each $t>0$
%
\begin{equation}
\sup_NE\biggl[\int_0^t\psi
\bigl(V^N(s)\bigr)\,ds\biggr]<\infty.\label{occbnd0}
\end{equation}
Then $\{\Gamma_N\}$ is relatvely compact, and if $\Gamma_N\Rightarrow
\Gamma$, then for
$f_1,\ldots, f_m\in D_{\psi}$,
\begin{eqnarray*}
&&\biggl(\int_0^{\cdot}f_1
\bigl(V^N(s)\bigr)\,ds,\ldots,\int_0^{\cdot}f_m
\bigl(V^N(s)\bigr)\,d s\biggr)
\\
&&\qquad\Rightarrow\biggl(\int
_{\Bbb E}f_1(v)\Gamma\bigl(dv\times[0,\cdot]\bigr),
\ldots,\int_
{\Bbb E}f_m(v)\Gamma\bigl(dv\times[0,
\cdot]\bigr)\biggr)
\end{eqnarray*}
in $C_{{\Bbb R}^m}[0,\infty)$.
\end{lemma}

\begin{pf}
Relative compactness of $\{\Gamma_N\}$ follows from Lemma 1.3 of
\citet{Ku92}. Relative compactness in $C_{{\Bbb R}^m}[0,\infty)$ follows
from relative compactness of each component. To see
that for $f\in D_{\psi}$, the sequence $X^N=\break \int_0^{\cdot}f(V^N(s
))\,ds$ is relatively compact, it
is enough to approximate the sequence by sequences
known to be relatively compact. For $\varepsilon>0$, there exists
a compact $K_{\varepsilon}\subset{\Bbb E}$ and $C>0$, such that $|f
|\leq(C{\bf1}_{K_{\varepsilon}}+\varepsilon)\psi$.
Define $X^N_{\varepsilon}=\int_0^{\cdot}{\bf1}_{K_{\varepsilon}}(V^N(s
))f(V^N(s))\,ds$. Note that $X_{\varepsilon}^N$ is
Lipschitz with Lipschitz constant $\sup_{v\in K_{\varepsilon}}|f(v)|$,
so $\{X_{
\varepsilon}^N\}$ is
relatively compact.
For $\delta>0$,
\[
\sup_NP\Bigl\{\sup_{s\leq t}\bigl|X^N(s)-X^N_{\varepsilon}(s)\bigr|
\geq\delta\Bigr\}\leq \frac{
\varepsilon}{\delta}\sup_NE\biggl[\int
_0^t\psi\bigl(V^N(s)\bigr)\,ds\biggr],
\]
and relative compactness of $\{X^N\}$ follows; see Problem
3.11.18 of \citet{EK86}.

Assuming that $\Gamma_N\Rightarrow\Gamma$, the convergence of $\int_
0^{\cdot}f(V^N(s))\,ds$ to
$\int_{{\Bbb E}\times[0,\cdot]}f(v)\times \Gamma(dv\times ds)$ follows by
the same type of
approximation.
\end{pf}

The final condition insures relative compactness of $V^N$.
%
\begin{condition}[(Tightness)]\label{tight}
If $\Bbb E$ is unbounded, there exists a locally bounded \mbox{$\psi
\dvtx {\Bbb E}\rightarrow[1,\infty
)$} satisfying
$\lim_{v\rightarrow\infty}\psi(v)=\infty$ such that for each $t>
0$,
%
\begin{equation}
\sup_NE\biggl[\int_0^t\psi
\bigl(V^N(s)\bigr)\,ds\biggr]<\infty\label{occbnd}
\end{equation}
and all of the following functions are in $D_{\psi}$: $\sup_N|F^N|$, $\sup_N|r_N(F^N-F)|$, $\sup_N|r_N(F-\bar{F}-A_NH_
N)|$,
$|G|$, $\sup_N|A_Nh|$ for $h\in{\cal D}(L_0)\cap B({\Bbb E}_0)$,
$\sup_
N|\frac{1}{r_{1,N}}A_Nh|$ for
$h\in{\cal D}(L_1)\cap B({\Bbb E}_0\times{\Bbb E}_1)$, and $\sup_
N|\frac{1}{r_{2,N}}A_Nh|$ for
$h\in{\cal D}(L_2)\cap B({\Bbb E})$.
\end{condition}

Assuming the above conditions and defining
%
\begin{equation}
\bar{G}(v_0)=\int\!\!\int G(v_0,v_1,v_2)
\mu_{v_0,v_1} (dv_2)\mu_{v_0}(dv_1),\label{gbar}
\end{equation}
and similarly for $\bar{G}_0$ and $\bar{G}_1$, we have the following
functional central limit theorem.

\begin{theorem}\label{genclt}
Under the above conditions, suppose that\break  $\lim_{N\rightarrow\infty}
U^N(0)=U(0)$, that
$\bar{F}$ is continuously differentiable and that
the solution (necessarily unique) of (\ref{limdeq})
exists for all time.
Then for each $t>0$,
\[
\sup_{s\leq t}\bigl|V_0^N(s)-V_0(s)\bigr|
\Rightarrow0,
\]
$r_N(M^{N,1}-M^{N,2})\Rightarrow M$, where $M$ has Gaussian, mean-zero,
independent increments with
%
\begin{equation}
E\bigl[M(t)M^T(t)\bigr]=\int_0^t
\bar{G}\bigl(V_0(s)\bigr)\,ds,\label{covar}
\end{equation}
and $U^N\Rightarrow U$ satisfying
\[
U(t)=U(0)+M(t)+\int_0^t \bigl(\nabla\bar{F}
\bigl(V_0(s)\bigr)U(s)+\bar{G}_
0\bigl(V_0(s)
\bigr)+\bar{G}_1\bigl(V_0(s)\bigr) \bigr)\,ds.
\]
Assuming $\bar{G}=\sigma\sigma^T$, we can write
%
\begin{eqnarray}\label{itogs}
U(t)&=&U(0)+\int_0^t\sigma\bigl(V_0(s)
\bigr)\,dW(s)
\nonumber
\\[-8pt]
\\[-8pt]
\nonumber
&&{}+\int_0^t \bigl(\nabla \bar{F}
\bigl(V_0(s)\bigr)U(s)+\bar{G}_0\bigl(V_0(s)
\bigr)+\bar{G}_1\bigl(V_0(s)\bigr) \bigr)\,ds.
\end{eqnarray}
\end{theorem}

\begin{remark}
As noted above, the corresponding theorem for systems
with two time-scales is obtained by assuming ${\Bbb E}_2$
consists of a single point so
$L_2f\equiv0$.
\end{remark}

\begin{pf*}{Proof of Theorem \ref{genclt}}
Let $\Gamma_N$ be the occupation measure defined as in~(\ref{occdef}).
Then by Lemma \ref{relcomp}, $\{\Gamma_N\}$ is relatively compact.
Assume, for simplicity that $\Gamma_N\Rightarrow\Gamma$. We will
show that $
\Gamma$
is uniquely determined.

Condition \ref{cnd1}, equation (\ref{jmp1}) and the martingale
central limit theorem, Theorem~\ref{mgclt}, imply
$M^{N,1}\Rightarrow0$, and Lemma \ref{relcomp} then implies
$V^N_0\Rightarrow V_0^{\infty}$, where
%
\begin{equation}
V^{\infty}_0(t)=V_0(0)+\int_{{\Bbb E}\times[0,t]}
F(v)\Gamma(dv\times ds).\label{v0lim}
\end{equation}

Condition \ref{tight}, the definition of $L_2$, and Lemma
\ref{relcomp}
imply
\begin{eqnarray*}
&&\frac{1}{r^{2,N}}\biggl(h\bigl(V^N(t)\bigr)-h\bigl(V^N(0)
\bigr)-\int_0^tA_Nh
\bigl(V^N(s)\bigr)\,ds\biggr)\\
&&\qquad\Rightarrow \int_{{\Bbb E}\times[0,t]}L_2h(v)
\Gamma(dv\times ds)
\end{eqnarray*}
for every $h\in C^{\infty}_c({\Bbb E})$. The uniform integrability implied
by (\ref{occbnd}) implies that the limit is a continuous
martingale with sample paths of finite variation and
hence is identically zero. Condition \ref{qstat} then
implies [see Example 2.3 of \citet{Ku92}] that $\Gamma$ can be written
\[
\Gamma(dv\times ds)=\mu_{v_0,v_1}(dv_2)\Gamma^{0,1}(dv_0
\times dv_1\times ds).
\]
A similar argument gives
\begin{eqnarray*}
0&=&\int_{{\Bbb E}\times[0,t]}L_1h(v)\Gamma(dv\times ds)\\
&=&\int
_{
{\Bbb E}_0\times{\Bbb E}_1\times[0,t]}\bar{L}_1h(v_0,v_1)
\Gamma^{
0,1}(dv_0\times dv_1\times ds),
\end{eqnarray*}
which implies
\[
\Gamma^{0,1}(dv_0\times dv_1\times ds)=
\mu_{v_0}(dv_1)\Gamma^0( dv_0\times
ds).
\]
But the convergence of $V_0^N$ to $V_0^{\infty}$ implies
$\Gamma^0(dv_0\times ds)=\delta_{V_0^{\infty}(s)}(dv_0)\,ds$.\vadjust{\goodbreak}

Now (\ref{v0lim}) can be rewritten
%
\begin{equation}
V^{\infty}_0(t)=V_0(0)+\int_0^t
\bar{F}\bigl(V_0^{\infty}(s)\bigr)\,ds,\label{v0lim2}
\end{equation}
and it follows that $V_0^{\infty}=V_0$.

Similarly, (\ref{cov1}) now becomes
\[
(r_N)^2\bigl[V_0^N-H_N
\circ V^N\bigr]_t\Rightarrow\int_0^t
\bar{G}\bigl(V_0(s) \bigr)\,ds,
\]
and it follows that $r_N(M^{N,1}-M^{N,2})\Rightarrow M$ as desired.

Finally, the uniform integrability implied by (\ref{occbnd})
and Condition \ref{tight} allows interchange of limits and
integrals in the expansion of $U^N$ given in (\ref{uexp}),
and the convergence of $U^N$ to $U$ follows.
\end{pf*}

\section{Diffusion approximation}\label{sectdiff}
The functional central limit theorem, Theorem~\ref{genclt},
suggests approximating $V^N_0$ by $V_0+\frac{1}{r_N}U$. In turn, that
observation and~(\ref{itogs})
suggest approximating $V_0^N$ by a diffusion
process given by the It\^o equation
%
\begin{eqnarray}\label{dfapp}
D^N(t)&=&V^N_0(0)+\frac{1}{r_N}\int
_0^t\sigma\bigl(D^N(s)
\bigr)\,dW(s)
\nonumber
\\[-8pt]
\\[-8pt]
\nonumber
&&{} +\int_0^t \biggl(\bar F
\bigl(D^N(s)\bigr)+\frac{1}{r_N}\bar G_0
\bigl(D^N(s )\bigr)+\frac{1}{r_N}\bar G_1
\bigl(D^N(s)\bigr) \biggr)\,ds.
\end{eqnarray}

The approximation
\[
V_0^N\approx\hat{D}^N\equiv V_0+
\frac{1}{r_N}U
\]
is, of course, justified
by Theorem \ref{genclt}. Justification for the
approximation $V_0^N\approx D^N$ is less clear, since $D^N$ is not
produced as a limit. Noting, however, that
\begin{eqnarray*}
\hat{D}^N(t)&=&V_0^N(0)+\frac{1}{r_N}
\int_0^t\sigma\bigl(V_0(s)
\bigr)\,dW(s)
\\
&& {}+\int_0^t \biggl(\bar F
\bigl(V_0(s)\bigr)+\frac{1}{r_N}\nabla\bar F\bigl(V_
0(s)
\bigr)U(s)+\frac{1}{r_N}\bar G_0\bigl(V_0(s)\bigr)\\
&&\hspace*{151pt}\qquad{}+
\frac{1}{r_N}\bar G_1\bigl(V_0(s)\bigr) \biggr)\,ds,
\end{eqnarray*}
assuming smoothness of $\bar F$, $\bar{G}_0$ and $\bar{G}_1$, we see that
$r_N^2(D^N-\hat{D}^N)$ converges to $\hat{U}$ satisfying
\begin{eqnarray*}
\hat{U}(t)&=&\int_0^t\nabla\sigma
\bigl(V_0(s)\bigr)U(s)\,dW(s)
\\
&&{} +\int_0^t \biggl(\nabla\bar{F}
\bigl(V_0(s)\bigr)\hat{U}(s)+\frac{1}2 U^T(s)
\partial^2\bar{F}\bigl(V_0(s)\bigr)U(s)
\\
&&\hspace*{65pt}{} +\bigl(\nabla\bar{G}_0\bigl(V_0(s)\bigr)+\nabla\bar{
G}_1\bigl(V_0(s)\bigr)\bigr)U(s) \biggr)\,ds,
\end{eqnarray*}
and since the central limit theorem demonstrates that
the fluctuations of $V^N$ are of order $O(r_N^{-1})$, we see that
the difference between the two approximations $D^N$ and
$\hat{D}^N$ is negligible compared to these fluctuations.


\section{Markov chain models for chemical reactions}\label{sectmcmcr}

A reaction network is a chemical system involving
multiple reactions and chemical species. The kind of
stochastic model for a network that we will consider
treats the system as a continuous time Markov chain
whose state $X$ is a vector giving the number of
molecules $X_i$ of each species of type $i\in\I$ present. Each
reaction is modeled as a possible transition for the
state. The model for the $k$th reaction, for each $k\in\K$,
is determined by a vector of inputs $\nu_k$ specifying the
numbers of molecules of each chemical species that are
consumed in the reaction, a vector of outputs $\nu'_k$
specifying the numbers of molecules of each species that
are produced in the reaction, and a function of the state
$\lambda_k(x)$ that gives the rate at which the reaction occurs as
a function of the state. Specifically, if the $k$th reaction
occurs at time $t$, the change in $X$ is a vector of integer values
$\zeta_k=\nu'_k-\nu_k$.

Let
$R_k(t)$ denote the number of times that the $k$th reaction
occurs by time $t$. Then $R_k$ is a counting process with
intensity $\lambda_k(X(t))$ (called the \textit{propensity} in the chemical
literature) and can be written as
\[
R_k(t)=Y_k\biggl(\int_0^t
\lambda_k\bigl(X(s)\bigr)\,ds\biggr),
\]
where the $Y_k$ are independent unit Poisson processes.
The state of the system at time $t$ can be written as
\[
X(t)=X(0)+\sum_k\zeta_kR_k(t)=X(0)+
\sum_k\zeta_kY_k
\biggl(\int_0^t\lambda_
k\bigl(X(s)
\bigr)\,ds\biggr).
\]

In the stochastic version of the law of mass action, the
rate function is proportional to the number of ways of
selecting the molecules that are consumed in the
reaction, that is,
\[
\lambda_k(x)=\kappa'_k\prod
_i\nu_{ik}!\prod_i\pmatrix{{x_i}
\cr{\nu_{
ik}}}=\kappa'_k\prod
_ix_i(x_i-1)\cdots(x_i-
\nu_{ik}+1).
\]
Of course, physically, $|\nu_k|=\sum_i\nu_{ik}$ is usually assumed to
be less than or equal to two, but that does not play a
significant role in the analysis that follows.

A reaction network may exhibit behavior on multiple
scales due to the fact that some species may be present
in much greater abundance than others, and the rate
functions may vary over several orders of magnitude.
Following \citet{KK10}, we embed the model of interest in
a sequence of models indexed by a scaling parameter $N$.
The model of interest corresponds to a particular value
of the scaling parameter $N_0$. For each species
$i\in\I=\{1,\ldots,s\}$, we specify a parameter $\alpha_i\geq
0$ and
normalize the number of molecules by $N_0^{\alpha_i}$ defining
$Z^{N_0}_i(t)=N_0^{-\alpha_i}X_i(t)$ so that it is of $O(1)$. For each
reaction $k\in\K$, we specify another parameter $\beta_k$ and
normalize the reaction rate constant as $\kappa'_k=\kappa_kN_0^{\beta_
k}$ so
that $\kappa_k$ is of $O(1)$. One can observe this model on different
time scales as well,
by replacing $t$ by $tN_0^\gamma$, for some $\gamma\in{\Bbb R}$.
The model then becomes a Markov chain
on ${\Bbb E}^{N_0}=N_0^{-\alpha_1}{\Bbb Z}_+\times\cdots
\times N_0^{-\alpha_{s}}{\Bbb Z}_+$ which, when $N=N_0$, evolves
according to
\[
Z^N_i(t)=Z^N_i(0)+\sum
_kN^{-\alpha_i}\zeta_{ik}Y_k
\biggl(\int_0^tN^{\nu_
k\cdot\alpha+\beta_k+\gamma}
\lambda_k^N\bigl(Z^N(s)\bigr)\,ds\biggr)
\]
with
\[
\lambda_k^N(z)=\kappa_k\prod
_iz_i\bigl(z_i-N^{-\alpha_i}
\bigr)\cdots\bigl(z_i -(\nu_{ik}-1)N^{-\alpha_i}\bigr).
\]
If for some $i$, $\alpha_i>0$ and $\nu_{ik}>1$, then $\lambda_k^N$
varies with $
N$
but converges as $N\rightarrow\infty$. To simplify notation, we will
write $\lambda_k(z)$ rather than $\lambda_k^N$, but one should check that
the $N$-dependence is indeed negligible in the analysis that
we do. Defining $\Lambda_N=\operatorname{diag}(N^{-\alpha_1},\ldots
,N^{-\alpha_s})$, so $Z^N=\Lambda_NX$, let
\[
A_Nf(z)=\sum_kN^{\rho_k}
\lambda_k(z) \bigl(f(z+\Lambda_N\zeta_k)-f(z)
\bigr),
\]
where $\rho_k=\nu_k\cdot\alpha+\beta_k+\gamma$.
Since the change of time variable from $t$ to $tN^{\gamma}$ is
equivalent to scaling the generator by a factor of $N^{\gamma}$, we
initially take $\gamma$ to be zero. We subsequently consider the
behaviour of $Z^N$ on different time-scales $Z^N(\cdot N^\gamma)$.

To be precise regarding the domain of $A_N$, note that
because the jumps of $Z^N$ are uniformly bounded, if we
define $\tau_r^N=\inf\{t\dvtx |Z^N(t)|\geq r\}$, then for every continuous
function $f$,
\[
f\bigl(Z^N\bigl(t\wedge\tau_r^N\bigr)
\bigr)-f\bigl(Z^N(0)\bigr)-\int_0^{t\wedge\tau_r^N}A_Nf
\bigl( Z^N(s)\bigr)\,ds
\]
is a martingale.

For notational simplicity, assume that the $\alpha_i$ satisfy
$0\leq\alpha_1\leq\cdots\leq\alpha_{s}$, and let $d_{\circ}\geq
0$ satisfy $\alpha_i=0$, $i\leq d_{\circ}$
and $\alpha_i>0$, $i>d_{\circ}$.

%


%

To apply the results of Section~\ref{sectclt}, we identify
$r_N,r_{1,N},r_{2,N}$ from the reaction network and the
parameters $\{\alpha_i\},\{\beta_k\}$ as follows. Let
\[
m_2=\max\{\rho_k-\alpha_i\dvtx
\zeta_{ik}\neq0\},
\]
and define $r_{2,N}=N^{m_2}$. Then there exists a linear
operator $L_2$ such that for each compact $K\subset{\Bbb R}^{s}$,
\[
\lim_{N\rightarrow\infty}\sup_{z\in K\cap{\Bbb E}^N}
\biggl|\frac{1}{r_{2,N}}A_Nh(z)-L_2h(z)\biggr|=0,\qquad
h\in{\cal D}(L_2)=C^1\bigl({\Bbb R}^{
s}
\bigr).
\]
Depending on the relationship between $\rho_k$ and $\alpha_i$ for
$\zeta_{ik}\neq0$ and the time-scale parameter $\gamma$, the
limiting operator $L_2$ is either the
generator for a Markov chain, a differential operator, or
a combination of the two, which would be the generator
for a piecewise deterministic Markov process (PDMP).
We classify the reactions by defining
\[
\K_{2,\circ}=\{k\in\K\dvtx \rho_k=m_2\}
\]
and
\[
\K_{2,\bullet}=\{k\in\K\dvtx \rho_k-\alpha_i=m_2
\mbox{ for some }i \mbox{ with } \alpha_i>0, \zeta_{ik}
\neq0\}.
\]
For each $k\in\K_{2,\circ}\cup\K_{2,\bullet}$
define
\[
\zeta_{2,k}=\lim_{N\rightarrow\infty}N^{\rho_k-m_2}
\Lambda_N\zeta _k\in\mathbb{Z}.
\]
Note that throughout the paper $\zeta_{2,k}$ will denote the limiting
reaction vector, not to be confused with the single matrix entry $\zeta_{ik}$.
Then, for $h\in C^1({\Bbb R}^{s})$
%
\begin{equation}
L_2h(z)=\sum_{k\in\K_{2,\circ}}\lambda_k(z)
\bigl(h(z +\zeta_{2,k})-h(z) \bigr)+\sum_{k\in\K_{2,\bullet}}
\lambda _k(z)\nabla h(z)\cdot\zeta_{2,k}.\label{L2}
\end{equation}
Note that, although $\lambda_k(z)$ depends on all species types, the
dynamics defined by $L_2$ makes changes only due to reactions $\K
_{2,\circ}\cup\K_{2,
\bullet}$. In other words,
only the subnetwork defined by reactions $\K_{2,\circ}\cup\K
_{2,\bullet}$ is
relevant on the time-scale corresponding to $\gamma=-m_2$.
If $\K_{2,\bullet}$ is empty, the process corresponding to $L_2$ is a
Markov chain, and if $\K_{2,\circ}$ is empty, the process is just
the solution of an ordinary differential equation. If both
are nonempty, the process is piecewise deterministic in
the sense of \citet{Dav93}.

The process corresponding to $L_2$ can be obtained as the
solution of
\[
V_2(t)=V_2(0)+\sum_{k\in\K_{2,\circ}}
\zeta_{2,k}Y_k\biggl(\int_
0^t
\lambda_k\bigl(V_2(s)\bigr)\,ds\biggr)+\sum
_{k\in\K_{2,\bullet}}\zeta_{2,k}\int_0^t
\lambda_k\bigl(V_2(s)\bigr)\,ds,
\]
and assuming that $V_2$ does not hit infinity in finite
time, $Z^N(\cdot N^{-m_2})\Rightarrow V_2$.

The central limit theorem in Section~\ref{sectclt} assumes that the
state space is a product space and that the fast process ``averages
out'' one component. The state space on which functions in the domain
of $L_1$ in Condition \ref{multconv} are defined is such that every
function on it is contained in the kernal of $L_2$. In order to
separate the state space in this way, we need to identify the
combinations of species variables whose change on the fastest
time-scale $\gamma=-m_2$ is less than $O(1)$. This can be done with a
change of basis of the original state space as follows.

Let $S_{\boldsymbol\K}$ be a matrix whose columns are $\zeta_k, k\in
{\boldsymbol\K}$ for some subset ${\boldsymbol\K}\subset\K$. Then
$S_{\boldsymbol\K}$ is the \textit{stoichiometric matrix} associated
with the reaction subnetwork ${\boldsymbol\K}$. For the species types
whose behavior is discrete, $S_{\boldsymbol\K}$ gives the possible
jumps, while for the species whose behavior evolves continuously,
$S_{\boldsymbol\K}$ determines the possible paths. We will let
$\mathcal{R}(S_{\boldsymbol\K})=\operatorname{span}\{\zeta_k, k\in
{\boldsymbol\K}\}\subset{\Bbb R}^{s}$ denote the range of
$S_{\boldsymbol\K}$, called the \textit{stoichiometric subspace} of the
chemical reaction subnetwork ${\boldsymbol\K}$, and we will let
\[
\mathcal{N}\bigl(S_{\boldsymbol\K}^{T}\bigr)=\biggl\{\theta\in{\Bbb
R}^{s}\dvtx \sum_{i\in\I}
\theta_
i\zeta_{ik}=0\ \forall k\in{\boldsymbol\K}\biggr\}
\]
denote the null space of $S_{\boldsymbol\K}^{T}$ which is the
othogonal complement of ${\cal R}(S_{\boldsymbol\K})$. For each
initial value $z_0$ of the reaction system, $z_0+\mathcal
{R}(S_{\boldsymbol\K})$ defines the \textit{stoichiometric compatibility
class} of the system. Then both stochastically and deterministically
evolving components of the system must remain in the stoichiometric
compatibility class for all time $t>0$. The linear combinations of the
species $\theta\cdot X$ for $\theta\in\mathcal{
N}(S_{\boldsymbol\K}^{T})$ are conserved quantities; that is, they
are constant along the trajectories of the evolution of the reaction
subnetwork~${\boldsymbol\K}$.

On the time scale $\gamma=-m_2$, the fast subnetwork determined by
$L_2$ has the stoichiometric matrix $S_2$ whose columns are $\{\zeta
_{2,k},k\in\K_{2,\circ}\cup\K_{2,\bullet}\}$. Define ${\cal
N}(S_2^{T})$ as above, and note that $\theta\cdot V_2$, $\theta\in
{\cal N}(S_2^{T})$, are conserved quantities for the fast subnetwork,
that is, $\theta\cdot V_2(t)$ does
not depend on $t$. 
Let $s_2$ denote the dimension of ${\cal R}(S_2)$, and $s'_1=s-s_2$ be the
dimension of ${\cal N}(S_2^{T})$.
We now replace the natural state space of the process by ${\cal
N}(S_2^{T})\times{\cal R}(S_2)$, mapping the original processes onto
this product space by the orthogonal projection $\Pi_{{\cal
N}(S_2^{T})}\times\Pi_{{\cal R}(S_2)}$, that is,
\[
\bigl(V_1^{\prime N}(t),V_2^N(t)\bigr)=
\bigl(\Pi_{{\cal N}(S_2^{T})}Z^N(t),\Pi_{{\cal
R}(S_2)}Z^N(t)
\bigr).
\]

Note that the original coordinates have
different underlying state spaces $N^{-\alpha_i}{\Bbb Z}$; however, the
change of basis will combine only those coordinates with
the same scaling parameter $\alpha_i$.
To
see that this is the case, note that by the definition of
$\zeta_{2,k}$, $\zeta_{2,ik}\neq0$ and $\zeta_{2,jk}\neq
0$ implies $\alpha_i=\alpha_j$. It follows that
there is a basis $\theta_1,\ldots,\theta_{s'_1}$ for ${\cal
N}(S_2^{T})$ such that $\theta_{il}\neq0$
and $\theta_{jl}\neq0$ implies $\alpha_i=\alpha_j$, and we can take
this basis
to be orthonormal. We denote the common scaling parameter by $\alpha
_{\theta_
l}$. Let $\Theta_1$ be the matrix with rows $\theta_1^{T}, \ldots
,\theta_{s'_1}^{T}$ so that
$(\Theta_1 z)^T=(\theta_1\cdot z,\ldots,\theta_{s'_1}\cdot z)^{T}$
and the orthogonal projection is given by
\[
\Pi_{{\cal N}(S_2^{T})}=\Theta_1^{T}\Theta_1=
\sum_{l=1}^{s'_1}\theta_l
\theta_l^T.
\]

On the next time scale we only need to consider the dynamics of the
projection of the original process that is unaffected by the fast
subnetwork $V_1^{\prime N}=\Pi_{{\cal N}(S_2^{T})}\Lambda_NX$. Since $\Pi
_{{\cal N}(S_2^{T})}\Lambda_N=\Lambda_N\Pi_{{\cal N}(S_2^{T})}$, we have
\[
V_1^{\prime N}(t)=\Pi_{{\cal N}(S_2^{T})}Z^N( 0)+
\Lambda_N\sum_k\Pi_{{\cal N}(S_2^{T})}
\zeta_kY_k\biggl(N^{\rho_k} \int
_0^t\lambda_k^N
\bigl(Z^N(s)\bigr)\,ds\biggr).
\]
Note that $\Pi_{{\cal R}(S_2)}\zeta_k$ is not necessarily equal to
$\zeta_{2,k}$, nor is the other projection
$\Pi_{{\cal N}(S_2^{T})}\zeta_k=\zeta_k-\Pi_{{\cal R}(S_2)}\zeta
_k$ necessarily equal to $\zeta_k-\zeta_{2,k}$.
To identify the next time scale let
\[
m_1=\max\{\rho_k-\alpha_{\theta_l}\dvtx
\theta_l\cdot\zeta_k\neq0\} =\max\bigl\{
\rho_k-\alpha_i\dvtx (\Pi_{{\cal N}(S_2^{T})}
\zeta_k)_i\neq 0\bigr\},
\]
and define
$r_{1,N}=N^{m_1}$. Note that $m_1<m_2$.
Then there exists a linear operator $L_1$ such that for
each compact $K\subset{\Bbb R}^{s'_1}$,
%
\begin{equation}
\lim_{N\rightarrow\infty}\sup_{z\in K\cap
{\Bbb E}^N}\biggl|\frac{1}{r_{1,N}}A_Nh(z)-L_1h(z)\biggr|=0,\label{L1}\vadjust{\goodbreak}
\end{equation}
where $h\in{\cal D}(L_1)$ satisfies $h(z)=f(\theta_1\cdot z,\ldots,
\theta_{s_1'}\cdot z)$ for
$f\in C^1({\Bbb R}^{s'_1})$.
Define
\[
\K_{1,\circ}=\Bigl\{k\in\K\dvtx \rho_k=m_1,\max
_l|\theta_l\cdot\zeta_k|> 0
\Bigr\}
\]
and
\[
\K_{1,\bullet}=\{k\in\K\dvtx \rho_k-\alpha_{\theta_l}=m_1
\mbox{ for some } l\mbox{ with }\alpha_{\theta_l}>0,
\theta_l\cdot\zeta_k\neq0 \}.
\]
%
Let $\Lambda_N^{\Theta_1}=\operatorname{diag}(N^{-\alpha_{\theta
_1}},\ldots
,N^{-\alpha_{\theta_{s'_1}}})$, and for each $k\in\K_{1,\circ}\cup
\K_{1,\bullet}$ define
\[
\zeta_{1,k}^{\theta}=\lim_{N\rightarrow\infty}N^{\rho
_k-m_1}
\Lambda_
N^{\Theta_1}\Theta_1\zeta_k=
\lim_{N\rightarrow\infty}N^{\rho_k-m_1}\bigl(N^{
-\alpha_{\theta_1}}
\theta_1\cdot\zeta_k,\ldots,N^{-\alpha
_{\theta_{
s'_1}}}
\theta_{s'_1}\cdot\zeta_k\bigr)^T.
\]
%
Then for $h(z)=f(\Theta_1 z)$ with $f\in C^1({\Bbb R}^{s'_1})$
\[
L_1h(z)=\sum_{k\in\K_{1,\circ}}
\lambda_k(z) \bigl(f\bigl(\Theta_1 z+\zeta
_{1,k}^{\theta}\bigr)-f(\Theta_1 z)\bigr)+\sum
_{k\in\K_{1,\bullet}}\lambda _k(z)\nabla f(
\Theta_1 z)\cdot\zeta_{1,k}^{\theta}.
\]
If $V_1$ denotes the process corresponding to $L_1$, 
then assuming that $V_1$ does not hit infinity in finite
time, $V_1^{\prime N}(\cdot N^{-m_1})=\Pi_{{\cal N}
(S_2^{T})} Z^N(\cdot N^{-m_1})\Rightarrow V_1$.

To separate the state space in terms of the next time scale (if there
is one), define
\[
\zeta_{1,k}=\lim_{N\rightarrow\infty}N^{\rho_k-m_1}
\Lambda_
N\Pi_{{\cal N}(S_2^{T})}\zeta_k.
\]
In other words, $\zeta_{1,k}=\Theta_1^{T}\zeta_{1,k}^{\theta}$ is
embedded in the original space, and $\zeta_{1,k}^{\theta}=\Theta
_1\zeta_k^1$.
On the time scale $\gamma=-m_1$, the subnetwork determined by $L_1$
has the stoichiometric
matrix $S_1$ with columns $\{\zeta_{1,k} k\in\K_{1,\circ}\cup\K
_{1,\bullet}
\}$. Define the subspace ${\cal N}(S_1^{T})$ as before, and let $s_1$
denote the dimension of ${\cal R}(S_1)$ and $s_0=s'_1-s_1$ be the
dimension of ${\cal N}(S_1^{T})$.
As before we need to map the processes $V_1^{\prime N}$ onto this product
space by the orthogonal projection $\Pi_{{\cal N}(S_1^{T})}\times\Pi
_{{\cal R}(S_1)}$. Since $\zeta_{1,k}\in{\cal N}(S_2^{T})=\operatorname{span}(\theta_1,\ldots,\theta_{
s'_1})$, we can assume that the
$\theta_l$ are selected so that
\[
{\cal R}(S_1)=\operatorname{span}(\theta_{s_0+1},\ldots,
\theta _{s'_1})=\operatorname{span} (\zeta_{1,1},\ldots,
\zeta_{1,s_1}).
\]
Define
\[
\Pi_0=\sum_{l=1}^{s_0}
\theta_l\theta_l^{T}=\Pi_{{\cal
N}(S_1^{T})},\qquad
\Pi_1=\sum_{l
=s_0+1}^{s'_1}
\theta_l\theta_l^{T}=\Pi_{{\cal R}(S_1)}\quad
\mbox{and}\quad \Pi_2=\Pi_{{\cal R}(S_2)}.
\]
%

On the next time scale we need only consider the projection $\Pi_0Z^N$
of the original process which is unaffected by either of the faster
subnetworks. To identify the next time scale, let
\[
m_0=\max\{\rho_k-\alpha_{\theta_l}\dvtx
\theta_l\cdot\zeta_k\neq0,1 \leq l\leq s_0
\}=\max\bigl\{\rho_k-\alpha_i\dvtx (\Pi_0
\zeta_k )_i\neq0\bigr\},
\]
and define
$r_{0,N}=N^{m_0}$. Note that if $1\leq s_2, 1\leq s_1, 1 \leq s_0$
($s_0+s_1+s_2=s$), $m_0<m_1<m_
2$.
Without loss of generality, we can assume that time is
scaled so that $m_0=0$.
Then, there exists a linear operator $L_0$ such that for
each compact $K\subset{\Bbb R}^{s_0}$,
\[
\lim_{N\rightarrow\infty}\sup_{z\in K\cap{\Bbb E}^N}\bigl|A_
Nh(z)-L_0h(z)\bigr|=0,
\]
where $h\in{\cal D}(L_0)$ satisfies $h(z)=f(\theta_1\cdot z,\ldots
,\theta_{s_0}\cdot z)$ for
$f\in C^1({\Bbb R}^{s_0})$.
Define
\[
\K_{0,\circ}=\Bigl\{k\in\K\dvtx \rho_k=0,\max
_l|\theta_l\cdot\zeta_k|>0\Bigr
\}
\]
and
\[
\K_{0,\bullet}=\{k\in\K\dvtx \rho_k-\alpha_{\theta_l}=0
\mbox{ for some } l\mbox{ with }\alpha_{\theta_l}>0,
\theta_l\cdot\zeta_k\neq0,1\leq l\leq s_0
\}.
\]
As before, let $\Theta_0$ be the matrix with rows
$\theta_1^T,\ldots,\theta_{s_0}^T$, and let $\Lambda_N^{\Theta
_0}=\break \operatorname{diag}
(N^{-\alpha_{\theta_1}}, \ldots,N^{-\alpha_{\theta_{s_0}}})$, so
that $\Pi_0=\Pi_{{\cal N}(S_1^{T})}=\Theta_0^{T}\Theta_0$ and for
each $k\in\K_{0,\circ}\cup\K_{0,\bullet}$ define
\[
\zeta_k^{\theta,0}=\lim_{N\rightarrow\infty}N^{\rho_k}
\Lambda_
N^{\Theta_0}\Theta_0\zeta_k=\lim
_{N\rightarrow\infty}N^{\rho_k}\bigl(N^{
-\alpha_{\theta_1}}
\theta_1\cdot\zeta_k,\ldots,N^{-\alpha
_{\theta_{
s_0}}}
\theta_{s_0}\cdot\zeta_k\bigr)^T.
\]
%
For $h(z)=f(\Theta_0z)$ with $f\in C^1({\Bbb R}^{s_0})$
\[
L_0h(z)=\sum_{k\in\K_{0,\circ}}\lambda_k(z)
\bigl(f\bigl(\Theta_0z+ \zeta_k^{\theta,0}\bigr)-f(
\Theta_0z)\bigr)+\sum_{k\in\K_{0,\bullet}}
\lambda_k(z)\nabla f(\Theta_0z)\cdot
\zeta_k^{\theta,0}.
\]
%
To relate the above calculations to the results of Section~\ref{sectclt}, we assume that $\K_{0,\circ}=\varnothing$ so that
\[
L_0h(z)=\sum_{k\in\K_{0,\bullet}}\lambda_k(z)
\nabla f(\Theta_
0z)\cdot\zeta_k^{\theta,0}.
\]

Let $V^N=TZ^N\equiv(\Pi_0Z^N,\Pi_1Z^N,\Pi_2Z^N)$,
so $V^N=(V_0^N,V_1^N,V_2^N)\in\break  {\cal N}(S_1^T)\times{\cal
R}(S_1)\times{\cal R}(S_2)$, and note that $T$
is invertible so that the intensities can be written as
functions of $v\in{\cal N}(S_1^T)\times{\cal R}(S_1)\times{\cal
R}(S_2)$, that is, $\lambda_k(T^{-1}v)$.
Since $\Pi_0z=\sum_{l=1}^{s_0}(\theta_l\cdot z)\theta_l$ and $\Pi_
1z=\sum_{l=s_0+1}^{s'_1}(\theta_l\cdot z)\theta_l$,
the process $V_0^N=\Pi_0Z^N$ is the embedding of $\Theta_0Z^N$, and
similarly $(V_0^N,V_1^N)=(\Pi_0Z^N,\Pi_1^NZ^N)$ is just the embedding
of $\Theta_1 Z^N$. Let
${\Bbb E}_0$, ${\Bbb E}_1$, and ${\Bbb E}_2$ denote the limit of the
state spaces for
$V_0^N$, $V_1^N$ and $V_2^N$.

The function $F^N$ in (\ref{mg1}) is given by
%
\begin{equation}
F^N(v)=\sum_kN^{\rho_k}
\Lambda_N^{\Theta_0}\lambda _k\bigl(T^{-1}v
\bigr)\Theta_0 \zeta_k\label{fndef}
\end{equation}
and
\[
F(v)=\lim_{N\to\infty}F^N(v)=\sum
_{k\in\K_{0,\bullet
}}\lambda_k\bigl(T^{-1}v\bigr)
\zeta_
k^{\theta,0}.
\]

To satisfy Condition~\ref{qstat} we will
assume that $L_2$ is such that for each
$(v_0,v_1)\in{\Bbb E}_0\times{\Bbb E}_1$ there exists a unique conditional
equilibrium distribution $\mu_{v_0,v_1}(dv_2)\in\mathcal{P}({\Bbb E}_
2)$ for $L_2$.
Then
$\bar{L}_1h(v_0,v_1)=\int L_1h(v_0,v_1,u_2)\mu_{v_0,v_1}(du_2)$ is
%
\begin{eqnarray}\label{L1bar}
\bar{L}_1h(v_0,v_1)&=&\sum
_{k\in\K_{1,\circ}} \bar{\lambda}_k(v_0,v_1)
\bigl(f\bigl((v_0,v_1)+\zeta_{1,k}^{\theta}
\bigr)-f(v_0,v_1) \bigr)
\nonumber
\\[-8pt]
\\[-8pt]
\nonumber
&&{}+\sum
_{k\in\K_{1,\bullet}}\bar{\lambda}_k(v_0,v_1)
\nabla f(v_0,v_1)\cdot\zeta_{1,k}^{\theta},
\end{eqnarray}
where $\bar{\lambda}_k(v_0,v_1)=\int\lambda_k(T^{-1}(v_0,v_1,v_2)
)\mu_{v_0,v_1}(dv_2)$.
For Condition~\ref{qstat} to be met, we also need to
assume that for each $v_0\in{\Bbb E}_0$ there exists a unique
conditional equilibrium distribution
$\mu_{v_0}(dv_1)\in\mathcal{P}({\Bbb E}_1)$ for $\bar{L}_1$.

We further need to assume that there are functions
$h_1\in\mathcal{D}(L_1)\dvtx {\Bbb E}_0\times{\Bbb E}_1\mapsto{\Bbb R}^{
|{\Bbb E}_0|}$ and
$h_2,h_3\in\mathcal{D}(L_2)\dvtx {\Bbb E}\mapsto{\Bbb R}^{|{\Bbb E}_0
|}$ that solve the following
Poisson equations:
\[
\bar{L}_1h_1=\bar{F}_1-\bar{F},\qquad
L_2h_2=F-\bar{F}_1,\qquad L_2h_3=
\bar{L}_1h_1-L_1h_1,
\]
where
\[
\bar{F}_1(v_0,v_1)=\int
F(v_0,v_1,u_2)\mu_{v_0,v_1}(du_2),\qquad
\bar{F}(v_0)=\int\bar{F}_1(v_0,u_1)
\mu_{v_0}(du_1)
\]
in order for Condition~\ref{cnd7} to be met. We refer the reader to
\citet{GM96} for results on sufficient conditions for the existence of
solutions to a Poisson equation for a general class of Markov
processes. For the class of general piecewise deterministic processes
see also \citet{CD03}. For the examples considered in Section~\ref{sectexamp}, we were able to explicitly compute the desired functions. In general,
however, explicit computation may not be possible, so results that
ensure the existence of these functions may be useful.\looseness=1

%



We now need to identify $r_N$, which will be of the form
$r_N=N^p$, for some $0<p<m_1$. Assuming that
there is no cancellation among the terms in the sum in~(\ref{fndef}), for (\ref{err1}) to hold, we must have
%
\begin{equation}
p\leq\max\{\alpha_{\theta_l}-\rho_k\dvtx \theta_l
\cdot \zeta_k\neq0,\rho_k<\alpha_{\theta_l},1\leq
l\leq s_0\}.\label{pbnd1}
\end{equation}
Then
\[
\theta_l\cdot G_0(v)=\lim_{N\to\infty}r_N
\theta_l\cdot \bigl(F^N(v)-F(v) \bigr)=\sum
_{k\dvtx \alpha_{\theta_l}-\rho_k=p}\lambda_k\bigl(T^{
-1}v\bigr)
\theta_l\cdot\zeta_k
\]
and
\[
G_0(v)=\sum_{l=1}^{s_0}\sum
_{k:\alpha_{\theta_l}-\rho_k=p}\lambda_
k\bigl(T^{-1}v
\bigr)\theta_l\cdot\zeta_k\theta_l.
\]

Now let $H^N=r_{1,N}^{-1}h_1+r_{2,N}^{-1}(h_2+h_3)$. To ensure that the
limit in
(\ref{err2}) exists, with reference to the definition of $L_2$,
we
must have
%
\begin{equation}
p\leq\min\{\alpha_i+m_2-\rho_k\dvtx
\zeta_{ik}\neq 0,\alpha_
i+m_2-
\rho_k>0\}\label{pbnd2}
\end{equation}
and
%
\begin{equation}
p\leq\min\{2\alpha_i+m_2-\rho_k\dvtx
\zeta_{ik}\neq0, \alpha_i>0\},\label{pbnd3}
\end{equation}
and with reference to the definition of $L_1$, we must have
%
\begin{equation}
p\leq\min\{\alpha_{\theta_l}+m_1-\rho_k\dvtx
\theta _l\cdot \zeta_k\neq0,\alpha_{\theta_l}+m_1-
\rho_k>0\}\label{pbnd4}
\end{equation}
and
%
\begin{equation}
p\leq\min\{2\alpha_{\theta_l}+m_1-\rho_k\dvtx
\theta_l \cdot\zeta_k\neq0,\alpha_{\theta_l}>0
\}.\label{pbnd5}
\end{equation}
Note that (\ref{pbnd1}) implies the minimum in (\ref{pbnd4}) and
(\ref{pbnd5}) only needs to be taken over $s_0+1\leq l\leq s'_1$.

Assuming that $h_1$, $h_2$ and $h_3$ are sufficiently smooth,
these assumptions insure that there
exists $G_1\dvtx {\Bbb E}\mapsto{\Bbb R}^{|{\Bbb E}_0|}$
\begin{eqnarray*}
G_1(v)&=&\lim_{N\to\infty} \biggl(r_N \biggl(
\frac
{A^N}{r_{2,N}}-L_2 \biggr) (h_2+h_3)+r_N
\biggl(\frac{A^N}{r_{1,N}}-L_
1 \biggr)h_1
\biggr)\\
&=&G_{12}( v)+G_{11}(v).
\end{eqnarray*}
To identify $G_{12}$,
define
\begin{eqnarray*}
\tilde{\zeta}_{2,k}&=&\lim_{N\rightarrow\infty}N^p
\bigl(N^{\rho
_k-m_2}\Lambda_
N\zeta_k-
\zeta_{2,k}\bigr),
\\
\tilde{\xi}_{2,kij}&=&\lim_{N\rightarrow\infty}N^{p+\rho
_k-m_2-\alpha_i-\alpha_
j}
\zeta_{ik}\zeta_{kj}
\end{eqnarray*}
and
\begin{eqnarray*}
\K_{2,\circ}^{p}&=&\{k\in\K\dvtx \theta_l\cdot
\zeta_k\neq0 \mbox{ for some $l$ with }\alpha_{\theta_l}=0,
m_2-\rho_k=p\},
\\
\K_{2,\bullet}^{p}&=&\{k\in\K\setminus\K_{2,\bullet}\dvtx
\theta _l\cdot\zeta_k\neq0 \mbox{ for some $l$ with }
\alpha_{\theta_l}>0, m_2-\rho_k+
\alpha_{\theta_l}=p\}.
\end{eqnarray*}
Then setting $h(z)=h_2(Tz)+h_3(Tz)$, $G_{12}(v)=H_{12}(T^{-1}v)$, where
\begin{eqnarray*}
H_{12}(z)&=&\sum_{k\in\K_{2,\circ}}
\lambda_k(z)\nabla h(z+\zeta _{2,k})\cdot\tilde{
\zeta}_{2,k}+\sum_{k\in\K_{2,\bullet}\cup\K
_{2,\bullet}^{p}}
\lambda_k(z)\nabla h(z)\cdot\tilde{\zeta}_{2,k}
\\
&&{} +\sum_{k\in\K_{2,\bullet}}\lambda_k(z)
\frac{1}2\sum_{ij}\partial_{z_i}\,
\partial_{z_j}h(z)\tilde{\xi}_{2,kij}+\sum
_{k\in\K_{2,\circ}^{p}}\lambda_k(z) \bigl(h(z+\tilde{
\zeta}_{2,k})-h(z)\bigr).
\end{eqnarray*}
Similarly, to identify $G_{11}$, define
\begin{eqnarray*}
\tilde{\zeta}_{1,k}^{\theta}&=&\lim_{N\rightarrow\infty
}N^p
\bigl(N^{\rho_k
-m_1}\Lambda^{\Theta_1}_N\Theta_1
\zeta_k-\zeta^\theta_{1,k}\bigr),
\\
\tilde{\xi}_{1,kll'}^{\theta}&=&\lim_{N\rightarrow\infty
}N^{p+\rho_k-m_1-\alpha_{
\theta_l}-\alpha_{\theta_{l'}}}
\zeta^{\theta}_{1,kl}\zeta^{\theta}_{1,
kl'}
\end{eqnarray*}
and
\begin{eqnarray*}
\K_{1,\circ}^{p}&=&\{k\in\K\dvtx \theta_l\cdot
\zeta_k\neq0 \mbox{ for some $l$ with }\alpha_{\theta_l}=0,
m_1-\rho_k=p\},
\\
\K_{1,\bullet}^{p}&=&\bigl\{k\in\K\setminus\K_{\bullet}^1
\dvtx \theta _l\cdot\zeta_k\neq0 \mbox{ for some $l$
with }\alpha_{\theta_l}>0, m_1-\rho_k+
\alpha_{\theta_l}=p\bigr\}
\end{eqnarray*}
Then $G_{11}(v)=H_{11}(T^{-1}v)$, where
\begin{eqnarray*}
H_{11}(z)&=&\sum_{k\in\K_{1,\circ}}
\lambda_k(z)\nabla h_1\bigl( \Theta_1 z+
\zeta^\theta_{1,k}\bigr)\cdot\tilde{\zeta}_{1,k}^{\theta
}\\
&&{}+
\sum_{k\in
\K_{1,\bullet}\cup\K_{1,\bullet}^{p}}\lambda_k(z)\nabla
h_1(\Theta_1 z)\cdot\tilde{\zeta}_{1,k}^{\theta}
\\
&&{} +\sum_{k\in\K_{1,\bullet}}\lambda_k(z)
\frac{1}2\sum_{ij}\partial_
l\,
\partial_{l'}h_1(\Theta_1 z)\tilde{
\xi}^{\theta}_{1,kll'}\\
&&{}+\sum_{k\in\K_{1,\circ}^{p}}
\lambda_k(z) \bigl(h_1\bigl(\Theta_1 z+
\tilde{\zeta }_{1,k}^{\theta}\bigr)-h_1(
\Theta_1 z)\bigr).
\end{eqnarray*}

We now need to identify $G\dvtx {\Bbb E}\to{\Bbb M}^{|{\Bbb E}_0|\times
|{\Bbb E}_0|}$ satisfying
(\ref{cov1}) in Condition \ref{hmgcnv}. Let
\[
R_k^N(t)=Y_k\biggl(N^{\rho_k}\int
_0^t\lambda_k\bigl(Z^N(s)
\bigr)\,ds\biggr)
\]
and $\tilde{H}_N(V^N)=\Theta_0Z^N-H_N(V^N)=V_0^N-H_N(V^N)$. Then
denoting $z^{\otimes
2}=zz^T$,
\begin{eqnarray*}
&& N^{2p}\bigl[\tilde{H}_N\bigl(V^N\bigr)
\bigr]_t\\
&&\qquad=\sum_kN^{2p}\int
_0^t\bigl(\tilde{H}_N
\bigl(V^N(s -)+T\Lambda_N\zeta_k\bigr)-
\tilde{H}_N\bigl(V^N(s-)\bigr)\bigr)^{\otimes2}\,dR_k^N(s)
,
\end{eqnarray*}
which is asymptotic to
\[
\sum_kN^{2p+\rho_k}\int_0^t
\bigl(\tilde{H}_N\bigl(V^N(s-)+T\Lambda_N
\zeta_
k\bigr)-\tilde{H}_N\bigl(V^N(s-)
\bigr)\bigr)^{\otimes2}\lambda_k\bigl(Z^N(s)\bigr)\,ds.
\]
Taking the limit as $N\to\infty$ and integrating with respect to $\mu
_{v_0,v_1}(dv_2)$ and $\mu_{v_0}(dv_1)$ then gives the value of $\bar{G}$.


\section{Examples}\label{sectexamp}

We now apply the central limit theorem to several
examples of chemical reaction networks with multiple
scales.

\subsection{Three species viral model}

Ball et~al. (\citeyear{BKPR06}) considered asymptotics for a
model of an intracellular viral infection originally given
in \citet{SYSY02} and studied further in \citet{HR02}.
The model includes three time-varying species, the viral
template,\vadjust{\goodbreak} the viral genome and the viral structural
protein, involved in six reactions,\vspace*{6pt}
\begin{center}
\begin{tabular}{lccc}
(1)& $T+\operatorname{stuff}$&$\stackrel{\kappa_1}{\rightharpoonup}$&$
T+G$,\\
(2)& $G$&$\stackrel{\kappa_2}{\rightharpoonup}$&$T$,\\
(3)& $T+\operatorname{stuff}$&$\stackrel{\kappa_3}{\rightharpoonup}$&$
T+S$,\\
(4)& $T$&$\stackrel{\kappa_4}{\rightharpoonup}$&$\varnothing$,\\
(5)& $S$&$\stackrel{\kappa_5}{\rightharpoonup}$&$\varnothing$,\\
(6)& $G+S$&$\stackrel{\kappa_6}{\rightharpoonup}$&$V$,
\end{tabular}\vspace*{6pt}
\end{center}
whose reaction rates (\textit{propensities}) are of
mass-action kinetics form $\lambda_k(x)=\kappa_k'\prod_ix_i^{\nu_{
ki}}$ with
constants\vspace*{6pt}
\begin{center}
\begin{tabular}{lcc}
$\kappa'_1$ & 1&  $1$,\\
$\kappa'_2$ & 0.025&  $2.5N_0^{-2/3}$,\\
$\kappa'_3$ & 1000&  $N_0$,\\
$\kappa'_4$ & 0.25&  $0.25$,\\
$\kappa'_5$ & $2$&  $2$,\\
$\kappa'_6$ & $7.5\times10^{-6}$&  $0.75N_0^{-5/3}$,
\end{tabular}\vspace*{6pt}
\end{center}
here expressed in terms of $N_0=1000$.

We denote $T$, $G$, $S$ as species 1, 2 and 3, respectively,
and let $X_i(t)$ denote the number of molecules of species $i$
in the system at time $t$. The stochastic model is
\begin{eqnarray*}
X_1(t)&=&X_1(0)+Y_2\biggl(\int
_0^t0.025X_2(s)\,ds
\biggr)\\
&&{}-Y_4\biggl(\int_0^t0.25X_1(s)
\,ds\biggr),
\\
X_2(t)&=&X_2(0)+Y_1\biggl(\int
_0^tX_1(s)\,ds\biggr)\\
&&{}-Y_2
\biggl(\int_0^t0.025X_2(s)\,ds\biggr)\\
&&{}-
Y_6\biggl(\int_0^t7.5
\cdot10^{-6}X_2(s)X_3(s)\,ds\biggr),
\\
X_3(t)&=&X_3(0)+Y_3\biggl(\int
_0^t1000X_1(s)\,ds
\biggr)\\
&&{}-Y_5\biggl(\int_0^t2X_3(s)\,ds
\biggr)\\
&&{}- Y_6\biggl(\int_0^t7.5
\cdot10^{-6}X_2(s)X_3(s)\,ds\biggr).
\end{eqnarray*}
We take
\[
\alpha_1=0,\qquad \alpha_2=2/3, \qquad\alpha_3=1.
\]
The scaling of the rate constants gives

\begin{center}
\begin{tabular}{cccc}
$k$ & $\kappa_k$ & $\beta_k$ & $\rho_k$ \\
$1$ & $1$ & $0$ & $0$\\
$2$ & $2.5$ & $-2/3$ & $0$\\
$3$ & $1$ & $1$ & $1$\\
$4$ & $0.25$ & $0$ & $0$\\
$5$ & $2$ & $0$& $1$\\
$6$ & $0.75$ & $-5/3$ & $0$.
\end{tabular}
\end{center}

Changing time $t\rightarrow N^{2/3}t$, the normalized system becomes
\begin{eqnarray*}
Z^N_1(t)&=&Z^N_1(0)+Y_2
\biggl(\int_0^tN^{2/3}2.5Z^N_2(s)\,ds
\biggr)-Y_4\biggl(\int_0^
tN^{2/3}0.25Z^N_1(s)\,ds
\biggr),
\\
Z^N_2(t)&=&Z^N_2(0)+N^{-2/3}Y_1
\biggl(\int_0^tN^{2/3}Z^N_1(s)\,ds
\biggr)\\
&&{}-N^{-2/
3}Y_2\biggl(\int_0^tN^{2/3}2.5Z^N_2(s)\,ds
\biggr)
\\
&& -N^{-2/3}Y_6\biggl(\int_0^tN^{2/3}0.75Z^N_2(s)Z^N_
3(s)\,ds
\biggr),
\\
Z^N_3(t)&=&Z^N_3(0)+N^{-1}Y_3
\biggl(\int_0^tN^{5/3}Z^N_1(s)\,ds
\biggr)-N^{-1}Y_
5\biggl(\int_0^tN^{5/3}2Z^N_3(s)\,ds
\biggr)
\\
&& -N^{-1}Y_6\biggl(\int_0^tN^{2/3}0.75Z^N_2(s)Z^N_3(
s)\,ds\biggr).
\end{eqnarray*}
We assume that the initial value for $Z^N_2$ is chosen to
satisfy $Z_2(0)=\break  \lim_{N\to\infty}Z^N_2(0)\in(0,\infty)$.

In this model, there are only two time-scales, so we set
\[
m_1=\max\{\rho_k-\alpha_i\dvtx
\zeta_{ik}\neq0\}=\max\bigl\{\tfrac{2}3-0,\tfrac{2}3-
\tfrac{2}3,\tfrac{5}3-1,\tfrac{2}3-1\bigr\}=
\tfrac{2}3,
\]
and we have
$r_{1,N}=N^{2/3}$. We have $\zeta_{1,1}=0, \zeta_{1,2}=e_1, \zeta
_{1,3}=e_3,  \zeta_{1,4}=-e_1,\break  \zeta_{1,5}=-e_3, \zeta_{1,6}=0$. The
operator $L_1=\lim_{N\rightarrow\infty}N^{
-2/3}A_N$ is given by
\begin{eqnarray*}
L_1h(z)&=&\lambda_2(z) \bigl(h(z+e_1)-h(z)
\bigr)+\lambda_4(z) \bigl(h( z-e_1)-h(z) \bigr)\\
&&{}+ \bigl(
\lambda_3(z)-\lambda_5(z) \bigr)\del_{z_3}h(z
)
\end{eqnarray*}
and note that for smooth $h$,
%
\begin{equation}
N^{-2/3}A_Nh=L_1h+O\bigl(N^{-2/3}
\bigr).\label{l1err}
\end{equation}

Functions $h\in\operatorname{ker}(L_1)$ are functions of the coordinate $
z_2$ only, ${\Bbb E}_1={\cal R}(S_1)=\operatorname{span}\{e_1,e_3\}$ and
${\Bbb E}_0={\cal N}(S_1^{T})=\operatorname{span}\{e_2\}$.
Taking $h\in{\cal D}(L_0)=C^1({\Bbb E}_0)$,
\[
L_0h(z)=\lim_{N\rightarrow\infty}A_Nh(z)= \bigl(
\lambda_1(z)-\lambda_
2(z)-\lambda_6(z)
\bigr)\del_{z_2}h(z_2).
\]

Setting $V_0^N=Z_2^N$ and $V_1^N=(Z_1^N,Z_3^N)$, the compensator for
$V_0^N$ is
\[
F^N(z)=\lambda_1(z)-\lambda_2(z)-
\lambda_6(z),
\]
so
$F(z)=F^N(z)$ and $G_0(z)\equiv0$ in Condition \ref{cnd1}.

The process corresponding to
$L_1$ is piecewise deterministic with $Z_1$ discrete and
$Z_3$ continuous. For fixed $z_2$, with reference to Condition
\ref{qstat},
the conditional equilibrium distribution
satisfies
%
\begin{eqnarray}\label{stateq1}
&&\int \biggl[2.5z_2 \bigl(g(z_1+1,z_3)-g(z_1,z_3)
\bigr)\nonumber\\
&&\hspace*{12pt}{}+0.25z_1 \bigl(g(z_
1-1,z_3)-g(z_1,z_3)
\bigr)
\\
&&\hspace*{63pt}{} +(z_1-2z_3)\frac{\del{g}}{\del{
z_3}}(z_1,z_3)
\biggr]\mu_{z_2}(dz_1,dz_3)=0.\nonumber
\end{eqnarray}
Note that the marginal for $Z_1$ is Poisson($10z_2$), so
\[
\int z_1\mu_{z_2}(dz_1,dz_3)=10z_2.
\]
Taking $g(z_1,z_3)=z_3$ in (\ref{stateq1}), we see
\[
\int z_3\mu_{z_2}(dz_1,dz_3)=5z_2.
\]
These calculations imply that
the averaged value for the drift $F$ is
\[
\bar{F}(z_2)=\int \bigl(\lambda_1(z)-
\lambda_2(z)-\lambda_6(z) \bigr)\mu_{z_2}(dz_1,dz_3)=7.5z_2-3.75z_2^2,
\]
with $\nabla\bar{F}(z_2)=7.5-7.5z_2$. For the current example, we
will see that $\bar{F}$ and $\bar{G}$ in~(\ref{gbar}) can be
obtained without
explicitly computing with $\mu_{z_2}$.

With reference to
(\ref{hnsol}),
we look for a
solution $h_1$ to the Poisson equation
%
\begin{equation}
L_1h_1(z)=(z_1-2.5z_2-0.75z_2z_3)-
\bigl(7.5z_2-3.75z_2^
2\bigr).\label{peq1}
\end{equation}
Trying $h_1$ of the form $h_1(z)=z_1u_1(z_2)+z_3u_3(z_2)$, we have
\[
L_1h_1(z)=u_1(z_2)
(2.5z_2-0.25z_1)+u_3(z_2)
(z_1-2z_3)
\]
and equating the factors multiplying $z_1$ and $z_3$, we get
$u_1(z_2)=1.5z_2-4$ and $u_3(z_2)=0.375z_2$. Thus
$h_1(z)=z_1(1.5z_2-4)+z_3(0.375z_2)$ and $H^N(z)=N^{-2/3}h_1(z)$.\vadjust{\goodbreak}

Since the solution of (\ref{peq1}) is exact and (as we
shall see) $r_N=N^{1/3}$, by (\ref{l1err}), we have $G_1=0$ in
Condition \ref{cnd7}. With reference to Condition
\ref{hmgcnv}, (\ref{jmp1}) and (\ref{hneg}) are immediate.

The only restriction that remains to determine $r_N$ is
the asymptotic behavior of the quadratic
variation of $Z^N_2-H^N(Z^N)=Z_2^N-N^{-2/3}h_1(Z^N)$. Direct
calculation shows that to get a nontrivial $G$ in
(\ref{cov1}) we must take $r_N=N^{1/3}$. We then have
\begin{eqnarray*}
&&N^{2/3} \bigl[Z^N_2-H^N
\bigl(Z^N\bigr) \bigr]_t\\
&&\qquad=\sum
_{
k=1}^6N^{-2/3}\int_0^t
\bigl(\zeta_{2k}+h_1\bigl(Z^N(s-)
\bigr)-h_1\bigl(Z^N(s-)+\Lambda_
N
\zeta_k\bigr)\bigr)^2\,dR_k^N(s)
\\
&&\qquad\approx\int_0^tZ_1^N(s)\,ds+
\int_0^t\bigl(-1-1.5Z_2^N(s)+4
\bigr)^22.5Z_2^N( s)\,ds
\\
&&\qquad\quad{} +\int_0^t\bigl(1.5Z_2^N(s)-4
\bigr)^2 0.25Z_1^N(s)\,ds+\int
_0^t0.75Z_2^N
(s)Z_3^N(s)\,ds,
\end{eqnarray*}
where we observe that jumps by $R_3^N$ and $R_5^N$ do not
contribute to the limit. Dividing the equation for $Z_1^N$ by
$N^{2/3}$, we observe that
\[
\int_0^tZ_1^N(s)\,ds
\approx\int_0^t10Z_2^N(s)\,ds.
\]
Similarly, dividing the equation for $Z_3^N$ by $N^{2/3}$ we see
that
\[
\int_0^tZ_3^N(s)\,ds
\approx\frac{1}2\int_0^tZ_1^N(s)\,ds
\approx\int_
0^t5Z_2^N(s)\,ds,
\]
which in turn implies
\[
\int_0^tZ_2^N(s)Z_3^N(s)\,ds
\approx\int_0^t5Z_2^N(s)^2\,ds.
\]
It follows that $\bar{G}(z_2)$ is
\begin{eqnarray*}
&&10z_2+(3-1.5z_2)^22.5z_2+(4-1.5z_2)^22.5z_2+3.75z_2^2
\\
&&\qquad=72.5z_2-48.75z_2^2+11.25z_2^3.
\end{eqnarray*}
Let $Z_2$ be the solution of
\[
Z_2(t)=Z_2(0)+\int_0^t
\bigl(7.5Z_2(s)-3.75Z_2^2(s) \bigr)\,ds
\]
and $U^N=N^{1/3}(Z^N_2-Z_2)$. Then
\[
\sup_{s\leq t}\bigl|Z^N_2(s)-Z_2(s)\bigr|
\Rightarrow0 \quad\mbox{and}\quad U^N\Rightarrow U,
\]
where, for $W$ a standard Brownian motion, $U$ satisfies
\begin{eqnarray*}
U(t)&=&U(0)+\int_0^t\sqrt{72.5Z_2(s)-48.75Z_2(s)^2+11.25Z_2(s)^3}
\,dW(s)\\
&&{}+\int_0^t \bigl(7.5-7.5Z_2(s)
\bigr)U(s)\,ds.
\end{eqnarray*}

The corresponding diffusion approximation is
\begin{eqnarray*}
D^N(t)&=&Z^N_2(0)\\
&&{}+N^{-1/3}\int
_0^t\sqrt{72.5D^N(s)-48.75D^N(s)^2+1
1.25D^N(s)^3}\,dW(s)
\\
&&{}+\int_0^t\bigl(7.5D^N(s)-3.75D^N(s)^2
\bigr)\,ds.
\end{eqnarray*}

\begin{figure}

\includegraphics{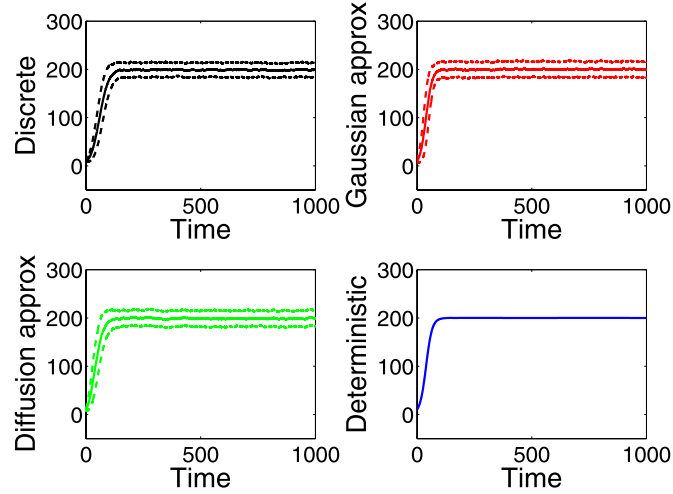}

\caption{Mean and standard deviation of the amount of genome in the
three species model
[$500$~simulations with parameters $N_0=1000$, $\gamma=2/3$,
$X_1(0)=0$, $X_2(0)=10$, $X_3(0)=0$].}
\label{three-cumul}
\end{figure}

\begin{figure}

\includegraphics{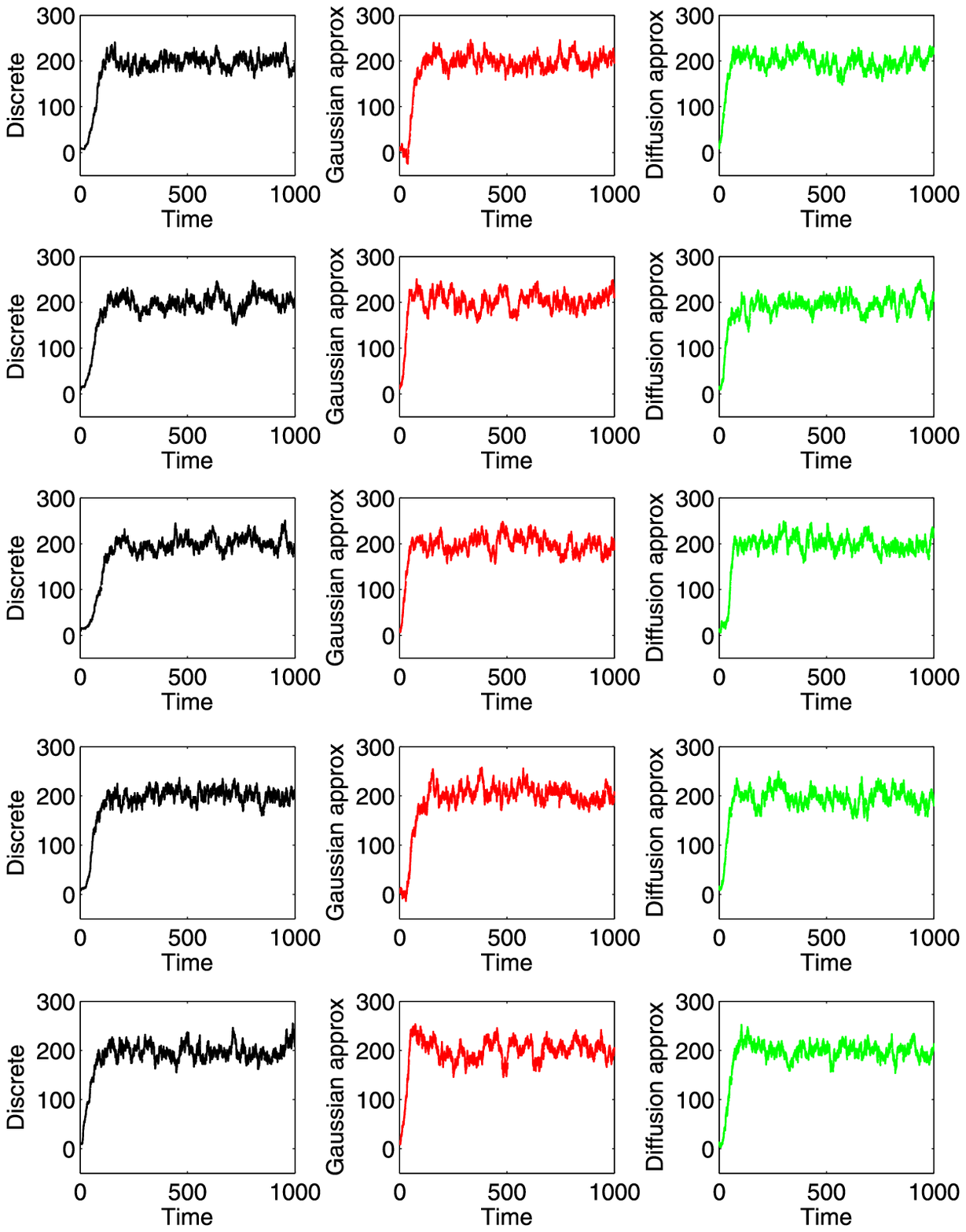}

\caption{Five trajectories of the amount of genome in the three
species model (same parameters as in Figure~\protect\ref{three-cumul}).}
\label{three-trajec}
\end{figure}

We compare simulations for the original value of the
amount of genome $X_2(\cdot)$ with the approximations given
by the Gaussian approximation\break
$N^{2/3}Z_2(\cdot N^{-2/3})+N^{1/3}U(\cdot N^{-2/3})$, and the diffusion
approximation\break  $N^{2/3}D^N(\cdot N^{-2/3})$. For comparison we also
give the deterministic value given by $N^{2/3}Z_2(\cdot N^{-2/3})$.
We use $N=1000$ and a time interval on the scale
$\gamma=2/3$. The initial values are set to
$X_1(0)=X_3(0)=0,X_2(0)=10$ and 500 realizations are
performed for each of the three stochastic processes.
Figure~\ref{three-cumul} shows the mean and one
standard deviation above and below the mean for each of
the three processes, and Figure~\ref{three-trajec} shows
five trajectories for the three processes.

For the diffusion process, these plots use only sample
paths that hit one \mbox{$(=100/N_0^{2/3})$} before they hit zero. For small
initial values, the diffusion
approximation does not give a good approximation of the
probability of hitting zero (and hence absorbing at zero),
before (e.g.) hitting one. Let
\[
\tau_Z^N=\inf\bigl\{t>0\dvtx Z^N_2(t)=0
\mbox{ or }Z_2^N(t)\geq1\bigr\}
\]
and
\[
\tau_D^N=\inf\bigl\{t>0\dvtx D^N(t)=0
\mbox{ or }D^N(t)\geq1\bigr\}.
\]
It is shown in \citet{BKPR06} that
\[
\lim_{N\rightarrow\infty}P\bigl\{Z^N\bigl(\tau_Z^N
\bigr)=0|Z^N(0)=N^{-2/3}k\bigr\}= 4^{-k}
\]
while a standard calculation for the diffusion process
gives
\[
\lim_{N\rightarrow\infty}P\bigl\{D^N\bigl(\tau_D^N
\bigr)=0|D^N(0)=N^{-2/3}k\bigr\}= e^{-({6}/{29})k}.
\]

\subsection{Michaelis--Menten enzyme model} A basic
model for an enzymatic reaction includes three
time-varying species, the substrate, the free enzyme
and the substrate-bound enzyme, involved in three
reactions:\vspace*{6pt}
\begin{center}
\begin{tabular}{lccc}
(1)& $S+E$&$\stackrel{\kappa'_1}{\rightharpoonup}$&$SE$,\\
(2)& $SE$&$\stackrel{\kappa_2'}{\rightharpoonup}$&$S+E$,\\
(3)& $SE$&$\stackrel{\kappa_3'}{\rightharpoonup}$&$P+E$,
\end{tabular}\vspace*{6pt}
\end{center}
with mass-action kinetics and with
rate constants such that $\kappa'_2,\kappa'_3\gg \kappa'_1$. To be precise,
let $\kappa'_2=\kappa_2N$, $\kappa'_3=\kappa_3N$, and $\kappa_1'=
\kappa_1$.

We denote $E$, $S$, $P$ as species 1, 2 and 3, respectively, and
let $X_i(t)$ be the number of molecules of species $i$ in the
system at time $t$. Note that the total number of
unbound and substrate-bound enzyme molecules is
conserved, and we let $M$ denote this amount. The
stochastic model is
\begin{eqnarray*}
X_1(t)&=&X_1(0)-Y_1\biggl(\int
_0^t\kappa'_1X_1(s)X_2(s)\,ds
\biggr)+Y_2\biggl(\int_0^t
\kappa'_2\bigl(M-X_1(s)\bigr)\,ds
\biggr)\\
&&{}+Y_3\biggl(\int_0^t
\kappa'_3\bigl(M-X_1(s)\bigr)\,ds\biggr),
\\
X_2(t)&=&X_2(0)-Y_1\biggl(\int
_0^t\kappa'_1X_1(s)X_2(s)\,ds
\biggr)+Y_2\biggl(\int_0^t
\kappa'_2\bigl(M-X_1(s)\bigr)\,ds\biggr),
\\
X_3(t)&=&X_3(0)+Y_3\biggl(\int
_0^t\kappa'_3
\bigl(M-X_1(s)\bigr)\,ds\biggr).
\end{eqnarray*}
If the initial amount of substrate is $O(N)\gg M$, then the
normalizations of the species abundances are given by
\[
\alpha_1=0,\qquad \alpha_2=1,\qquad \alpha_3=1,
\]
and the scaling exponents for the rate constants are
\[
\beta_1=0,\qquad \beta_2=1,\qquad \beta_3=1.
\]
The normalized system becomes
\begin{eqnarray*}
Z^N_1(t)&=&Z^N_1(0)-Y_1
\biggl(\int_0^tN\kappa_1Z^N_1(s)Z^N_2(s)\,ds
\biggr)+Y_2\biggl( \int_0^tN
\kappa_2\bigl(M-Z^N_1(s)\bigr)\,ds\biggr)
\\
&& +Y_3\biggl(\int_0^tN
\kappa_3\bigl(M-Z^N_1(s )\bigr)\,ds\biggr),
\\
Z^N_2(t)&=&Z^N_2(0)-N^{-1}Y_1
\biggl(\int_0^tN\kappa_1Z^N_1(s)Z^N_2(s)\,ds
\biggr)\\
&&{}+N^{-1}Y_2\biggl(\int_0^tN
\kappa_2\bigl(M-Z^N_1(s)\bigr)\,ds\biggr),
\\
Z^N_3(t)&=&Z^N_3(0)+N^{-1}Y_3
\biggl(\int_0^tN\kappa_3
\bigl(M-Z^N_1(s)\bigr)\,ds\biggr).
\end{eqnarray*}

Again, there are only two time-scales with the fast
time-scale
$m_1=1$ giving
$r_{1,N}=N$. Then $\zeta_{1,1}=-e_1, \zeta_{1,2}=\zeta_{1,3}=e_1$, and
the operator $L_1$ is given by
\[
L_1h(z)=\kappa_1z_1z_2
\bigl(h(z-e_1)-h(z) \bigr)+(\kappa_2+
\kappa_
3) (M-z_1) \bigl(h(z+e_1)-h(z)
\bigr),
\]
and for smooth $h$,
%
\begin{equation}
N^{-1}A_Nh=L_1h+O\bigl(N^{-1}
\bigr).\label{l1err2}
\end{equation}

Functions $h\in\operatorname{ker}(L_1)$ are functions of coordinates
$z_2$ and $z_3$ only. Thus ${\Bbb E}_1=\{z_1e_1\dvtx z_1=0,\ldots,M\}
\subset{\cal R}(S_1)$ and
${\Bbb E}_0={\cal N}(S_1^{T})=\{(z_2e_2,z_3e_3)\dvtx z_2,z_3\geq0\}$. For
$h\in{\cal D}
(L_0)=C^1({\Bbb E}_0)$,
\[
L_0h(z)=\bigl(\kappa_2(M-z_1)-
\kappa_1z_1z_2\bigr)\del_{z_2}h(z)+
\kappa_3 (M-z_1)\del_{z_3}h(z).
\]

Taking $V_0^N=(Z^N_2,Z^N_3)$, the compensator for $V_0^N$ in (\ref
{mg1}) is
\[
F^N(z)=\bigl(\kappa_2(M-z_1)-
\kappa_1z_1z_2,\kappa_3(M-z_1)
\bigr)^{T},
\]
so $F(z)=F^N(z)$ and $G_0(z)\equiv0$.

On the fast time-scale, the process whose generator is
$L_1$ is a Markov chain on ${\Bbb E}_1$ describing the dynamics of
an urn scheme with a total of $M$ molecules, and
for a fixed value of $z_2,z_3$, with transition rates $\kappa_1z_
2$ for
outflow and $\kappa_2+\kappa_3$ for inflow. Its stationary
distribution $\mu_{z_2,z_3}(z_1)$ is binomial$(M,p(z_2))$ for
\[
p(z_2)=\frac{\kappa_2+\kappa_3}{\kappa_2+\kappa_3+\kappa_1z_2},
\]
so $\int z_1\mu_{z_2,z_3}(dz_1)=Mp(z_2)$.

This observation
implies that the averaged value for the drift $F$ is
\begin{eqnarray*}
\bar{F}(z_2,z_3)&=&\biggl(-M\frac{\kappa_1\kappa_3z_2}{\kappa_2+\kappa_
3+\kappa_1z_2},M
\frac{\kappa_1\kappa_3z_2}{\kappa_2+\kappa
_3+\kappa_
1z_2}\biggr)^T\\
&=&-\kappa_3M
\bigl(1-p(z_2)\bigr)\pmatrix{1\cr{-1}}
\end{eqnarray*}
with
\[
\nabla\bar{F}=-M\frac{\kappa_1\kappa_3(\kappa_2+\kappa
_3)}{(\kappa_
2+\kappa_3+\kappa_1z_2)^2}\pmatrix{
 1&0
\cr
-1&0},
\]
and we need to solve the Poisson equation
\begin{eqnarray*}
L_1h_1(z)&=&\biggl(\kappa_2(M-z_1)-
\kappa_1z_1z_2+M\frac{\kappa_1\kappa_
3z_2}{\kappa_2+\kappa_3+\kappa_1z_2},\\
&&\hspace*{52pt}\kappa_3(M-z_1)-M\frac{\kappa_
1\kappa_3z_2}{\kappa_2+\kappa_3+\kappa_1z_2}\biggr)^T
\end{eqnarray*}
for $h_1$. Trying $h_1$ of the form $h_1(z)=(z_1u_1(z_2),z_1u_2(
z_2))^T$,
we have
\begin{eqnarray*}
L_1h_1(z)&=&\bigl(-\kappa_1z_1z_2u_1(z_2)+(
\kappa_2+\kappa_3) (M-z_1)u_
1(z_2),-
\kappa_1z_1z_2u_2(z_2)\\
&&\hspace*{140pt}{}+(
\kappa_2+\kappa_3) (M-z_1)u_2(z_2)
\bigr)^T,
\end{eqnarray*}
and equating terms with the same power of $z_1$, we get
$u_1(z_2)=(\kappa_1z_2+\kappa_2)/(\kappa_1z_2+\kappa_2+\kappa_3)$ and
$u_2(z_2)=\kappa_3/(\kappa_1z_2+\kappa_2+\kappa_3)$. Note that $
u_1(z_2)+ u_2(z_2)=1$.
Thus
\begin{eqnarray*}
h_1(z)&=&\biggl(\frac{z_1(\kappa_1z_2+\kappa_2)}{(\kappa_1z_2+\kappa_2
+\kappa_3)},\frac{z_1\kappa_3}{(\kappa_1z_2+\kappa_2+\kappa_3)}
\biggr)^
T\\
&=&z_1\bigl(u_1(z_2),1-u_1(z_2)
\bigr)^T,
\end{eqnarray*}
and $H^N(z)=N^{-1}h_1(z)$.

Examining the quadratic variation of $V_0^N-H^N\circ V^N$, we
see that $r_N$ must be $N^{1/2}$, and by (\ref{l1err2}), it
follows that $G_1=0$ in (\ref{err2}).

Finally, letting $z^{\otimes2}=zz^T$,
\begin{eqnarray*}
&&N \bigl[V_0^N-H^N\circ V^N
\bigr]_t\\
&&\qquad= N^{-1}\sum_{k=1}^3
\int_0^t \bigl(\Theta_0
\zeta_
k+h_1\bigl(Z^N(s-)
\bigr)-h_1\bigl(Z^N(s-)+\Lambda_N
\zeta_k\bigr) \bigr)^{\otimes
2}\,dR_k^N(s)
\\
&&\qquad\approx\int_0^t \biggl(-\pmatrix{1 \cr 0}+\pmatrix{{u_1\bigl(Z_2^N(s)\bigr)}
\vspace*{2pt}\cr{1-u_
1\bigl(Z_2^N(s)\bigr)}}
\biggr)^{\otimes2}\kappa_1Z_1^N(s)Z_2^N(s)\,ds
\\
&&\qquad\quad{} +\int_0^t \biggl(\pmatrix{1\cr 0}-\pmatrix{{u_1
\bigl(Z_2^N(s)\bigr)}\vspace*{2pt}\cr{1-u_
1
\bigl(Z_2^N(s)\bigr)}} \biggr)^{\otimes2}
\kappa_2\bigl(M-Z_1^N(s)\bigr)\,ds
\\
&&\qquad\quad{} +\int_0^t \biggl(\pmatrix{0\cr 1}-\pmatrix{{u_1
\bigl(Z_2^N(s)\bigr)}\vspace*{2pt}\cr{1-u_
1
\bigl(Z_2^N(s)\bigr)}} \biggr)^{\otimes2}
\kappa_3\bigl(M-Z_1^N(s)\bigr)\,ds
\\
&&\qquad\approx\int_0^t\pmatrix{ \bigl(1-u_1\bigl(Z_2^N(s)
\bigr)\bigr)^2&-\bigl(1-u_1\bigl(Z_2^N(s)
\bigr)\bigr)^2
\vspace*{2pt}\cr
-\bigl(1-u_1\bigl(Z_2^N(s)\bigr)
\bigr)^2&\bigl(1-u_1\bigl(Z_2^N(s)
\bigr)\bigr)^2 }
\kappa_1Z_1^N(s)Z_2^N(s)\,ds
\\
&&\qquad\quad{} +\int_0^t\pmatrix{ \bigl(1-u_1\bigl(Z_2^N(s)
\bigr)\bigr)^2&-\bigl(1-u_1\bigl(Z_2^N(s)
\bigr)\bigr)^2
\vspace*{2pt}\cr
-\bigl(1-u_1\bigl(Z_2^N(s)\bigr)
\bigr)^2&\bigl(1-u_1\bigl(Z_2^N(s)
\bigr)\bigr)^2 }
\kappa_2\bigl(M-Z_1^N(s)\bigr)\,ds
\\
&&\qquad\quad{} +\int_0^t\pmatrix{
 u_1\bigl(Z_2^N(s)
\bigr)^2&-u_1\bigl(Z_2^N(s)
\bigr)^2
\vspace*{2pt}\cr
-u_1\bigl(Z_2^N(s)\bigr)^2&u_1
\bigl(Z_2^N(s)\bigr)^2 }
\kappa_3\bigl(M-Z_1^N(s)
\bigr)\,ds,
\end{eqnarray*}
and averaging $Z_1^N$ gives
\begin{eqnarray*}
\lim_{N\rightarrow\infty}N \bigl[V_0^N-H_N
\circ V^N \bigr]_t&=&\int_
0^t
\bar{G}\bigl(Z(s)\bigr)\,ds\\
&=&\int_0^t\pmatrix
{ \bar{g}\bigl(Z_2(s)\bigr)&-\bar{g}
\bigl(Z_2(s)\bigr)
\vspace*{2pt}\cr
-\bar{g}\bigl(Z_2(s)\bigr)&\bar{g}\bigl(Z_2(s)\bigr)}\,ds,
\end{eqnarray*}
where $Z=(Z_2,Z_3)$ satisfies
\[
Z(t)=Z(0)+\int_0^tM\frac{\kappa_1\kappa_3Z_2(s)}{\kappa_2+\kappa_
3+\kappa_1Z_2(s)}\pmatrix{{-1}
\cr 1}\,ds
\]
and
\begin{eqnarray*}
\bar{g}(z_2)&=&M\bigl(1-u_1(z_2)
\bigr)^2\bigl(\kappa_1p(z_2)z_2+
\kappa_2\bigl(1-p(z_2) \bigr)\bigr)
\\
&&{}+Mu_1(z_2)^2\kappa_3
\bigl(1-p(z_2)\bigr).
\end{eqnarray*}
Let $U^N=N^{1/2}(Z^N_2-Z_2,Z^N_3-Z_3)^T$. Then
\[
\sup_{s\leq t}\bigl|\bigl(Z^N_2(s)-Z_2(s),Z_3^N(s)-Z_3(s)
\bigr)\bigr|\Rightarrow0 \quad\mbox{and}\quad U^N\Rightarrow U,
\]
where $U=(U_2,U_3)^T$ satisfies
\begin{eqnarray*}
U(t)&=&U(0)+\int_0^t\pmatrix {{-1}\cr 1}\sqrt{\bar{g}
\bigl(Z_2(s)\bigr)}\,dW(s)\\
&&{}+ \int_0^t
\frac{M\kappa_1\kappa_3(\kappa_2+\kappa_3)}{(\kappa
_2+\kappa_
3+\kappa_1Z_2(s))^2}U_2(s)\pmatrix{{-1}\cr 1}\,ds
\end{eqnarray*}
for $W$ a standard scalar Brownian motion.

The corresponding diffusion approximation is
\begin{eqnarray*}
\pmatrix{{D_2^N(t)}\cr{D_3^N(t)}}&=&\pmatrix{{Z_2^N(0)}
\cr{Z_3^N(0)}}+N^{
-1/2}\int
_0^t\pmatrix{{-1}\cr 1}\sqrt{\bar{g}
\bigl(D_2^N(s)\bigr)}\,dW(s)\\
&&{}+\int_
0^tM
\frac{\kappa_1\kappa_3D^N_2(s)}{\kappa_2+\kappa_3+\kappa_1D^
N_2(s)}\pmatrix{{-1}\cr 1}\,ds.
\end{eqnarray*}

\begin{figure}

\includegraphics{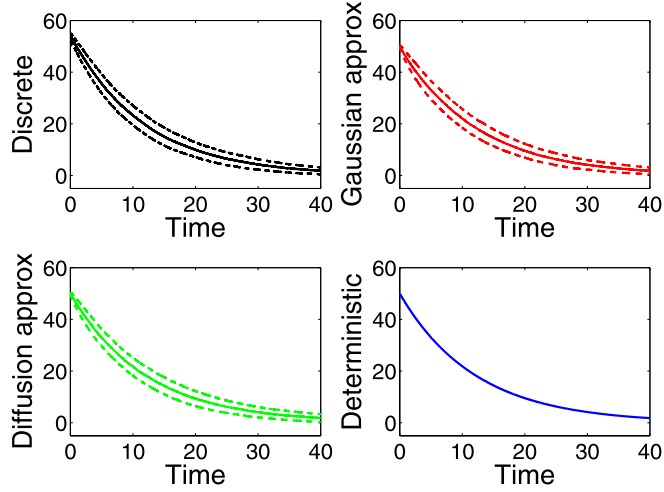}

\caption{Mean and standard deviation of the amount of substrate in the
Michaelis--Menten model
[$500$ simulations with parameters
$N_0=100$, $\gamma=0$, $X_1(0)=X_3(0)=0$, $X_2(0)=50$, $M=5$,
$\kappa_1'=0.1$, $\kappa_2'=500$, $\kappa_3'=100$].}
\label{mm-cumul}
\end{figure}

\begin{figure}

\includegraphics[scale=0.98]{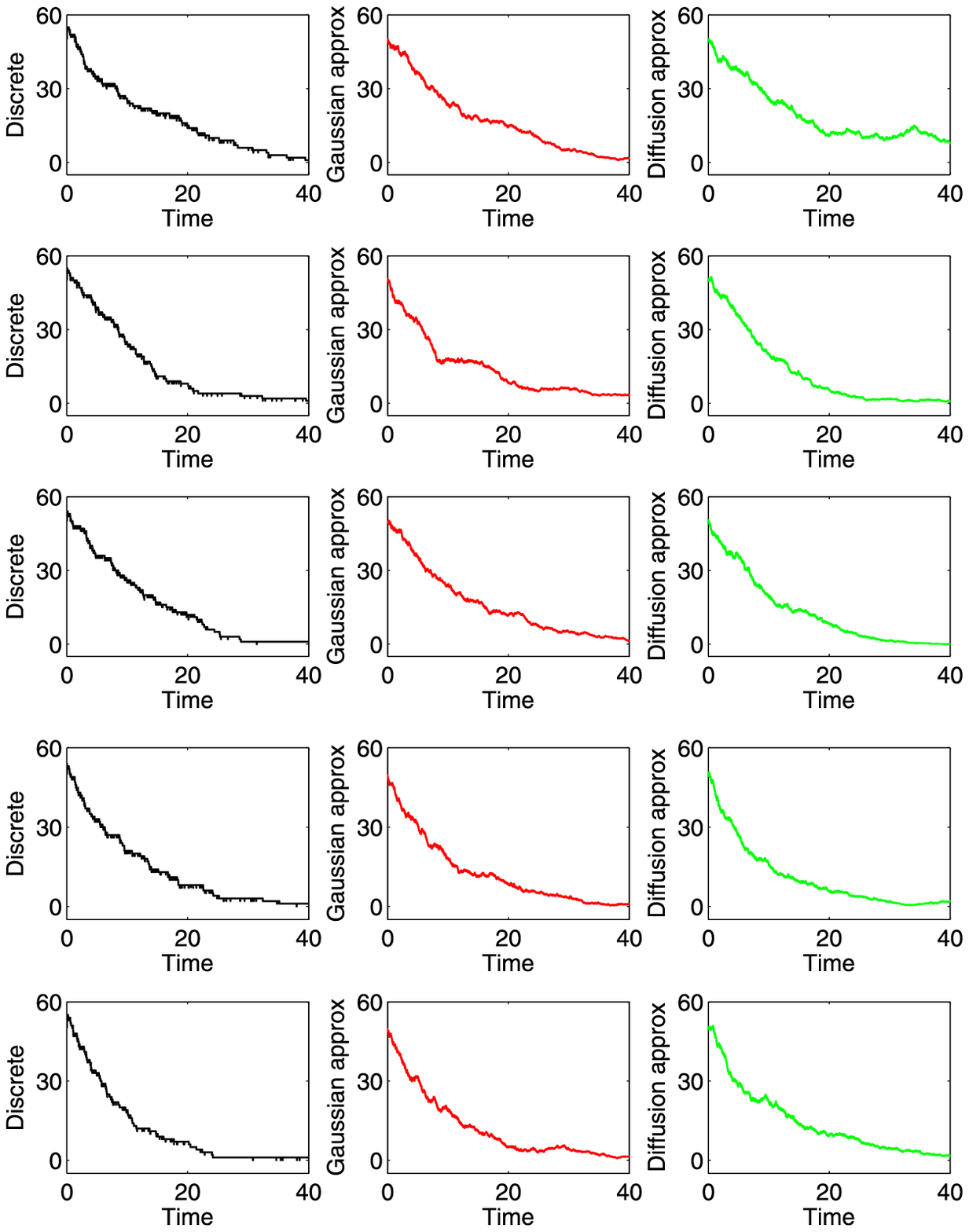}

\caption{Five trajectories of the amount of substrate in the
Michaelis--Menten model (parameters as in Figure~\protect\ref{mm-cumul}).}
\label{mm-trajec}\vspace*{-6pt}
\end{figure}

We compare simulations for 500 realizations of the
original model
with 500 realizations of the Gaussian approximation
$N_0Z_2(\cdot)+N_0^{1/2}U_2(\cdot),N_0Z_3(\cdot)+N_0^{1/2}U_3(
\cdot)$ and of the diffusion
approximation $N_0D^{N_0}_2(\cdot),N_0D^{N_0}_3(\cdot)$. For
comparison we also
give the deterministic value given by $N_0Z_2(\cdot),N_0Z_3(\cdot
)$. We
use $N_0=100$ and a time interval on the scale $\gamma=0$. The
initial values are set to $X_1(0)=X_3(0)=0,X_2(0)=50$ and
$M=5,\kappa_1'=0.1$, $\kappa_2'=500$ and $\kappa_3'=100$.
Figure~\ref{mm-cumul} shows the mean and one standard
deviation above and below the mean for each of the
three processes, and Figure~\ref{mm-trajec} shows five
trajectories for the three processes. In this example,
both Gaussian and diffusion approximations give good
approximations for the means and the standard
deviations of the pair of processes~$X_2(\cdot),X_3(\cdot)$.

\subsection{Another enzyme model} Another model for an
enzymatic reaction includes an additional form for the
enzyme which cannot bind to the substrate. There are
now four species, substrate, active enzyme,
enzyme-substrate complex and
inactive enzyme, involved in five reactions:
\begin{center}
\begin{tabular}{lccc}
(1)& $S+E$&$\stackrel{\kappa'_1}{\rightharpoonup}$&$SE$,\\
(2)& $SE$&$\stackrel{\kappa_2'}{\rightharpoonup}$&$S+E$,\\
(3)& $SE$&$\stackrel{\kappa_3'}{\rightharpoonup}$&$P+E$,
\end{tabular}
\end{center}
\begin{center}
\begin{tabular}{lccc}
(4)& $F$&$\stackrel{\kappa_4'}{\rightharpoonup}$&$E$,\\
(5)& $E$&$\stackrel{\kappa_5'}{\rightharpoonup}$&$F$,
\end{tabular}
\end{center}
with mass-action kinetics and
rate constants such that
$\kappa'_1=O(1)$, $\kappa'_2,\kappa'_3=O(N)$, $\kappa'_4,\kappa'_
5=O(N^2)$ so that
$\kappa_1'=\kappa_1$, $\kappa'_2=\kappa_2N$, $\kappa'_3=\kappa
_3N$, $\kappa'_
4=\kappa_4N^2$, $\kappa'_5=\kappa_5N^2$.\

We denote $E$, $S$, $F$ as species 1, 2 and 3, respectively,
and let $X_i(t)$ be the number of molecules of species $i$ in
the system at time $t$. The total number $M$ of active,
inactive and substrate-bound enzyme molecules is
conserved. The
stochastic model is
\begin{eqnarray*}
X_1(t)&=&X_1(0)-Y_1\biggl(\int
_0^t\kappa'_1X_1(s)X_2(s)\,ds
\biggr)\\
&&{}+Y_2\biggl(\int_0^t
\kappa'_2\bigl(M-X_1(s)-X_3(s)
\bigr)\,ds\biggr)
\\
&& {}+Y_3\biggl(\int_0^t
\kappa'_3\bigl(M-X_1(s)-X_3(s)
\bigr)\,ds\biggr)+Y_4\biggl(\int_0^t
\kappa'_
4X_3(s)\,ds\biggr)\\
&&{}-Y_5
\biggl(\int_0^t\kappa'_5X_1(s)\,ds
\biggr),
\\
X_2(t)&=&X_2(0)-Y_1\biggl(\int
_0^t\kappa'_1X_1(s)X_2(s)\,ds
\biggr)\\
&&{}+Y_2\biggl(\int_0^t
\kappa'_2\bigl(M-X_1(s)-X_3(s)
\bigr)\,ds\biggr),
\\
X_3(t)&=&X_3(0)-Y_4\biggl(\int
_0^t\kappa'_4X_3(s)\,ds
\biggr)+Y_5\biggl(\int_0^t
\kappa'_
5X_1(s)\,ds\biggr).
\end{eqnarray*}
If the initial amount of substrate is $O(N)\gg M$, then the
scaling exponents for the species abundances are
\[
\alpha_1=0,\qquad \alpha_2=1,\qquad \alpha_3=0,
\]
and the scaling exponents for the rate constants are
\[
\beta_1=0,\qquad \beta_2=1,\qquad \beta_3=1,\qquad
\beta_4=2, \qquad\beta_
5=2.
\]
The normalized system becomes
\begin{eqnarray*}
Z^N_1(t)&=&Z^N_1(0)-Y_1
\biggl(\int_0^tN\kappa_1Z^N_1(s)Z^N_2(s)\,ds
\biggr)\\
&&{}+Y_2\biggl( \int_0^tN
\kappa_2\bigl(M-Z^N_1(s)-Z^N_3(s)
\bigr)\,ds\biggr)
\\
&&{} +Y_3\biggl(\int_0^tN
\kappa_3\bigl(M-Z^N_1(s)-Z^N_3(s)
\bigr)\,ds\biggr)+Y_4\biggl(\int_0^tN^
2
\kappa_4Z^N_3(s)\,ds\biggr)\\
&&{}-Y_5
\biggl(\int_0^tN^2
\kappa_5Z^N_1(s)\,ds\biggr),
\\
Z^N_2(t)&=&Z^N_2(0)-N^{-1}Y_1
\biggl(\int_0^tN\kappa_1Z^N_1(s)Z^N_2(s)\,ds
\biggr)\\
&&{}+N^{-1}Y_2\biggl(\int_0^tN
\kappa_2\bigl(M-Z^N_1(s)-Z^N_3(s)
\bigr)\,ds\biggr),
\\
Z^N_3(t)&=&Z^N_3(0)-Y_4
\biggl(\int_0^tN^2
\kappa_4Z^N_3(s)\,ds\biggr)+Y_5
\biggl(\int_
0^tN^2
\kappa_5Z^N_1(s)\,ds\biggr).
\end{eqnarray*}

The fastest time-scale has $m_2=2$ and
$r_{2,N}=N^2$, with $\zeta_{2,4}=e_1-e_3,\zeta_{2,5}=-e_1+e_3$. The operator
$L_2$ is
\[
L_2h(z)=\kappa_4z_3 \bigl(h(z+e_1-e_3)-h(z)
\bigr)+\kappa_5z_1 \bigl( h(z-e_1+e_3)-h(z)
\bigr),
\]
with $\operatorname{ker}(L_2)$ consisting of functions of coordinates $
z_2$
and $z_1+z_3$ only. To simplify our calculations we make a
change of variables
to $(v_0,v_1,v_2)=(z_2,z_1+z_3,z_3)$, so
in this system of variables $\zeta_{2,4}=\tilde{e}_2,\zeta_{
2,5}=-\tilde{e}_2$ with
the operator $L_2$
\[
L_2h(v)=\kappa_4v_2 \bigl(h(v-
\tilde{e}_2)-h(z) \bigr)+\kappa_5(v_
1-v_2)
\bigl(h(v+\tilde{e}_2)-h(v) \bigr).
\]
Functions $h(v)\in\operatorname{ker}(L_2)$ are now functions of
$v_0,v_1$ only. Thus $\mathbb{E}_2={\cal R}(S_2)=\operatorname{span}
\{\tilde{e}_2\}$ and
$\mathbb{E}_1\times\mathbb{E}_0={\cal N}(S_2^T)=\operatorname{span}\{
\tilde{e}_1,\tilde{e}_0\}$.

The next time-scale has
$m_1=1$, $r_{1,N}=N$ and $\zeta_{1,1}=(0,-1), \zeta_{1,2}=\zeta
_{1,3}=(0,1)$. Also
\begin{eqnarray*}
L_1h(v)&=&
\kappa_1v_0(v_1-v_2)
\bigl(h\bigl((v_0,v_1-1)
-h(v_0,v_1)
\bigr) \bigr) 
\\
&&{}+(\kappa_2+\kappa_3)
(M-v_1) (h\bigl((v_0,v_1+1)
-h(v_0,v_1) \bigr)
\end{eqnarray*}
with $\operatorname{ker}(L_1)$ consisting of functions of $v_0$ only. Thus
$\mathbb{E}_1={\cal R}(S_1)=\operatorname{span}\{\tilde{e}_1\}$ and
$\mathbb{E}_1={\cal N}(S_1^T)=\operatorname{span}\{\tilde{e}_0\}$.

Finally, $L_0$ is
\[
L_0h(v)=-\kappa_1v_0(v_1-v_2)
\del_{v_0}h(v_0)+\kappa_2(M-v_1)
\del_{
v_0}h(v_0).
\]

The conditional stationary distribution $\mu_{v_0,v_1}(dv_2)$ of
Markov chain with generator $L_2$ is such that
$\rho_0(v_0,v_1)=\int v_2\mu_{v_0,v_1}(dv_2)=\frac{v_1\kappa
_5}{\kappa_4+\kappa_5}$,
thus 
\begin{eqnarray*}
\bar{L}_1h(v)&=&\kappa_1v_0
\frac{v_1\kappa_4}{\kappa_4+\kappa
_5} \bigl(h\bigl((v_0,v_1-1)
-h(v_0,v_1)\bigr) \bigr)\\
&&{}+(\kappa_2+
\kappa_3) (M-v_1) (h\bigl((v_0,v_1+1)
-h(v_0,v_1) \bigr),
\end{eqnarray*}
%
which has conditional stationary distribution $\mu_{v_0}(dv_1)$
such that
\begin{eqnarray*}
\rho_1(v_0)&=&\int v_1\mu_{v_0}(dv_1)=
\frac{M(\kappa_4+\kappa_5)
(\kappa_2+\kappa_3)}{\kappa_1\kappa_4v_0+(\kappa_4+\kappa
_5)(\kappa_
2+\kappa_3)},
\\
\rho_2(v_0)&=&
\int v_2
\mu_{v_0,v_1}(dv_2)\mu_{v_0}(dv_1) =
\frac{M
\kappa_5(\kappa_2+\kappa_3)}{\kappa_1\kappa_4v_0+(\kappa_4+\kappa_
5)(\kappa_2+\kappa_3)}.
\end{eqnarray*}
The compensator for the process $V^N_0$ is $F^N(v)=\kappa
_2(M-v_1)-\kappa_1v_0(v_1-v_2)=F(v)$, 
and averaging $F$ gives
$\bar{F}_1(v_0,v_1)=
\kappa_2(M-v_1)-\kappa_1v_0(v_1-\rho_0(v_0,v_1))$,
and
\begin{eqnarray*}
\bar{F}(v_0)
&=&\kappa_2\bigl(M-
\rho_1(v_0)\bigr)-\kappa_1v_0
\bigl(\rho_1(v_0)-\rho_2(v_0)
\bigr) 
\\
&=&-\frac{M\kappa_1\kappa_3\kappa_4v_0}{\kappa_1\kappa_4v_0+(\kappa_
4+\kappa_5)(\kappa_2+\kappa_3)},
\end{eqnarray*}
so
\[
\nabla\bar{F}(v_0)=-\frac{M\kappa_1\kappa_3\kappa_4(\kappa_4+
\kappa_5)(\kappa_2+\kappa_3)}{(\kappa_1\kappa_4v_0+(\kappa
_4+\kappa_
5)(\kappa_2+\kappa_3))^2}.
\]
Setting
\[
u_1(v_0)=\frac{\kappa_1\kappa_4v_0+\kappa_2(\kappa_4+\kappa
_5)}{\kappa_1\kappa_4v_0+(\kappa_2+\kappa_
3)(\kappa_4+\kappa_5)},\qquad u_2(v_0)=
\frac{\kappa_1v_0}{\kappa_4+\kappa_5},
\]
and
\[
u_3(v_0)=-\frac{\kappa_1v_0}{\kappa_4+\kappa_5}u_1(v_0)=-
\frac
{(\kappa_1\kappa_4v_0+\kappa_2(\kappa_4+\kappa_5))\kappa
_1v_0}{(\kappa_1\kappa_4v_0+(\kappa_2+\kappa_3)(\kappa_4+\kappa
_5))(\kappa_4+\kappa_5)},
\]
the solutions to the Poisson equations 
are given by functions
\[
h_1(v)=v_1u_1(v_0),\qquad
h_2(v)=-v_2u_2(v_0),\qquad
h_3(v)=-v_2u_3(v_0),
\]
and $H^N=\frac{1}Nh_1+\frac{1}{N^2}(h_2+h_3)$.

Let $r_N=N^{1/2}$
and observe that $\frac{1}{N^2}(h_2+h_3)$ makes a negligible
contribution to the quadratic variation. Consequently,
\begin{eqnarray*}
&&N \bigl[V_0^N-H^N\circ V^N
\bigr]_t\\
&&\qquad\approx\sum_{k=1}^5N^{-1}
\int_
0^t \bigl(\zeta_{k2}+h_1
\bigl(V^N(s-)\bigr)-h_1\bigl(V^N(s-)+T
\Lambda_N\zeta_k\bigr) \bigr)^2\,dR_k^N(s)
\\
&&\qquad\approx\int_0^t\bigl(-1+u_1
\bigl(V_0^N\bigr)\bigr)^2\kappa_1V_0^N
\bigl(V_1^N-V_2^N\bigr)\,ds
\\
&&\qquad\quad{} +\int_0^t\bigl(1-u_1
\bigl(V_0^N\bigr)\bigr)^2\kappa_2
\bigl(M-V_1^N\bigr)\,ds +\int_0^tu_1
\bigl(V_0^N\bigr)^2\kappa_3
\bigl(M-V_1^N\bigr)\,ds.
\end{eqnarray*}
Hence
\begin{eqnarray*}
G(v)&=&
 \bigl(\bigl(\kappa_3(\kappa_4+\kappa_5) \bigr)^2 \bigl(\kappa
_1v_0(v_1-v_2)+\kappa_2(M-v_1) \bigr)
\\
&&\hspace*{22pt}{}+ \bigl({\kappa_1\kappa_4v_0+\kappa_2(\kappa_4+\kappa_5)}
\bigr)^2 \bigl(\kappa_3(M-v_1) \bigr)\bigr)\\
&&{}/\bigl( \bigl(\kappa_1\kappa_4v_0+(\kappa
_2+\kappa_3)(\kappa_4+\kappa_5) \bigr)^2\bigr)
\end{eqnarray*}
and 
\[
\bar{G}(v_0)=\frac{M\kappa_1\kappa_3\kappa_4v_0 (\kappa
_3(\kappa_4+\kappa_5)^2(2\kappa_2+\kappa_3)+(\kappa_1\kappa
_4v_0+\kappa_2(\kappa_4+\kappa_5))^2 )}{(\kappa_1\kappa
_4v_0+(\kappa_2+\kappa_3)(\kappa_4+\kappa_5))^3}.
\]
If $V_0$ is the solution of
\[
V_0(t)=V_0(0)-\int_0^t
\frac{M\kappa_1\kappa_3\kappa
_4V_0(s)}{\kappa_
1\kappa_4V_0(s)+(\kappa_4+\kappa_5)(\kappa_2+\kappa_3)}\,ds,
\]
then, since $G_0=G_1\equiv0$, $U^N=N^{1/2}(V^N_0-V_0)\Rightarrow U$ where
\begin{eqnarray*}
U(t)&=&U(0)+\int_0^t\sqrt{\bar{G}
\bigl(V_0(s)\bigr)}\,dW_s\\
&&{}-\int_0^t
\frac{M
\kappa_1\kappa_3\kappa_4(\kappa_4+\kappa_5)(\kappa_2+\kappa_3)}{(
\kappa_1\kappa_4V_0(s)+(\kappa_4+\kappa_5)(\kappa_2+\kappa_3))^2} U(s)\,ds.
\end{eqnarray*}

The corresponding diffusion approximation is
\begin{eqnarray*}
D^N(t)&=&Z_2^N(0)+N^{-1/2}\int
_0^t\sqrt{\bar{G}\bigl(D^N(s)
\bigr)}\,dW(s)\\
&&{}-\int_
0^t\frac{M\kappa_1\kappa_3\kappa_4D^N(s)}{\kappa_1\kappa
_4D^N(s)+(\kappa_
4+\kappa_5)(\kappa_2+\kappa_3)}\,ds.
\end{eqnarray*}

\begin{figure}

\includegraphics{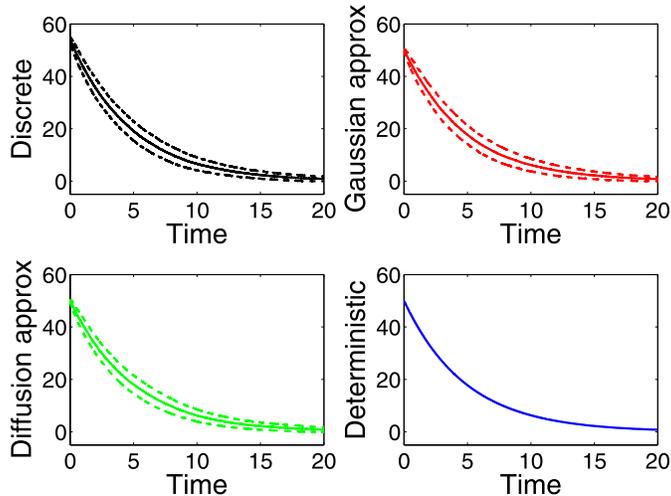}

\caption{Mean and standard deviations of the amount of substrate in
the three time-scale enzyme model
[$500$ simulations with $N_0=100$, $M=5$, $\gamma=0$, $X_1(0)=0$,
$X_2(0)=50$, $X_3(0)=0$,
$\kappa_1'=0.5$, $\kappa_2'=500$, $\kappa_3'=100$, $\kappa
_4'=\kappa_5'=5000$].}
\label{another-cumul}
\end{figure}

\begin{figure}

\includegraphics{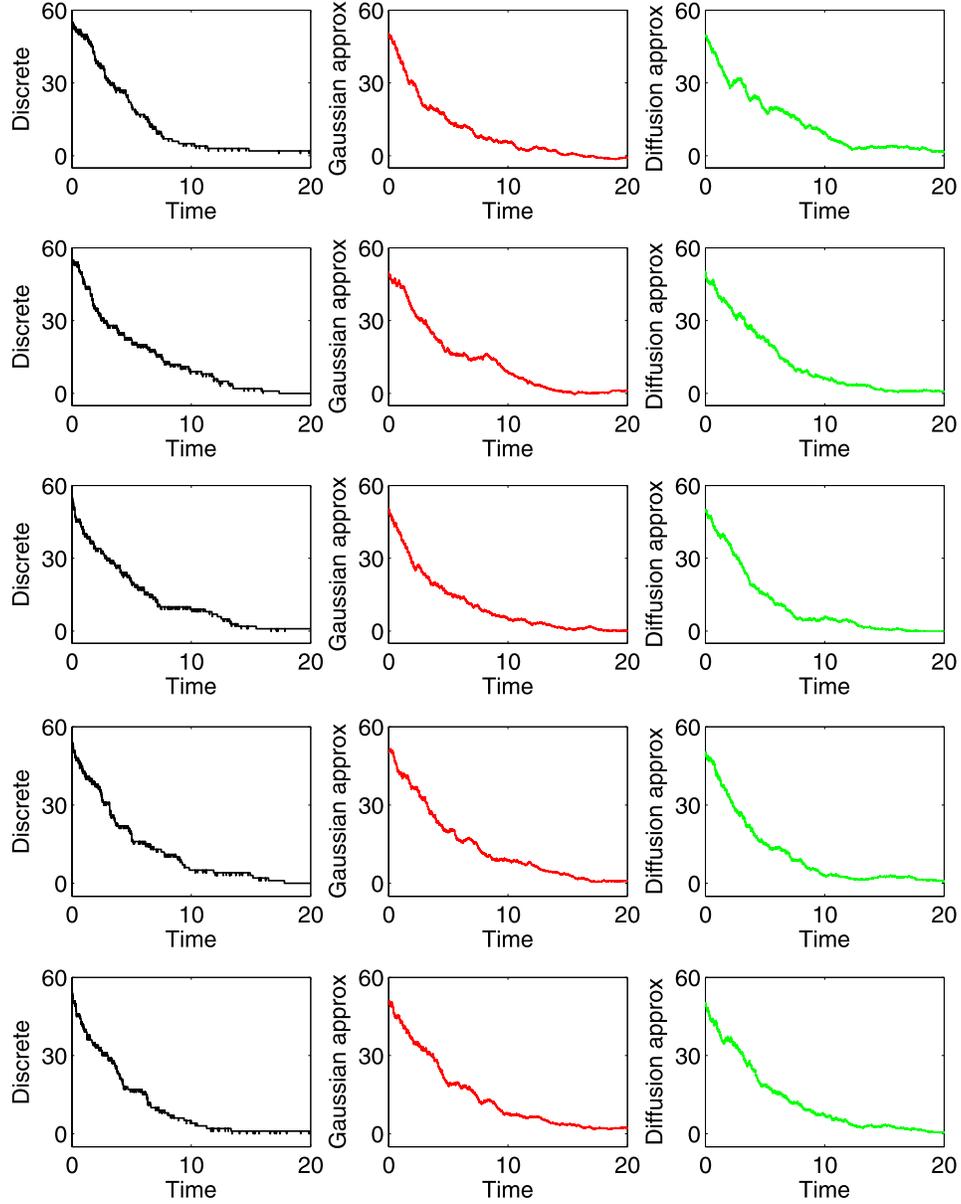}

\caption{Five trajectories for the amount of substrate in the three
time-scale enzyme model (same parameters as in Figure~\protect\ref
{another-cumul}).}
\label{another-trajec}
\end{figure}

Finally, we compare simulations for 500
realizations of the original model $X_2$
with 500 realizations of the
Gaussian approximation
$N_0V_0(\cdot)+N_0^{1/2}U(\cdot)$ and the diffusion approximation
$N_0D^{N_0}_2(\cdot)$. For comparison we also give the deterministic
value given by $N_0V_0(\cdot)$. We use $N_0=100$, a time interval
on the scale $\gamma=0$, and initial values are set to
$X_1(0)=X_3(0)=0,X_2(0)=50$ as in the previous example.
Here the additional parameters are set to $M=5$, $\kappa_1'=0.5$,
$\kappa_2'=500$, $\kappa_3'=100$ and $\kappa_4'=\kappa_5'=5000$.
Figure~\ref{another-cumul} shows the mean and one
standard deviation above and below the mean for each of
the three processes, and Figure~\ref{another-trajec} five
trajectories for the three processes. Again, both
Gaussian and diffusion approximations give a good
approximation for the mean and the standard deviation
from the mean of $X_2(\cdot)$.

\begin{appendix}

\setcounter{equation}{0}

\section*{Appendix}

\subsection{Martingale central limit theorem} Various
versions of the martingale central limit have been given
by \citet{McL74}, \citeauthor{Roo77} (\citeyear{Roo77,Roo80}), \citet{GH79} and
\citet{Reb80} among others. The following version is
from \citet{EK86}, Theorem 7.1.4.

\begin{theorem}\label{mgclt}
Let $\{M_n\}$ be a sequence of ${\Bbb R}^d$-valued martingales. Suppose
%
\begin{equation}
\lim_{n\rightarrow\infty}E\Bigl[\sup_{s\leq t}\bigl|M_n(s)-
M_n(s-)\bigr|\Bigr]=0\label{mgclt1}
\end{equation}
and
\[
\bigl[M^i_n,M^j_n
\bigr]_t\rightarrow c_{i,j}(t)
\]
for all $t\geq0$, where
$C=((c_{i,j}))$ is deterministic and continuous.
Then $M_n\Rightarrow M$, where $M$ is Gaussian with independent increments
and $E[M(t)M(t)^T]=C(t)$.
\end{theorem}

\begin{remark} Note that $C(t)-C(s)$ is nonnegative definite for
$t\geq s\geq0$. If $C$ is
absolutely continuous, then the derivative will also be nonnegative
definite and will have a nonnegative definite square
root. Suppose $\dot{C}(t)=\sigma(t)^2$ where $\sigma$ is symmetric.
Then $
M$ can
be written as
\[
M(t)=\int_0^t\sigma(s)\,dW(s),
\]
where $W$ is $d$-dimensional standard Brownian motion.
\end{remark}
\end{appendix}

%



\printaddresses

\end{document}